\documentclass[11pt]{amsart}

\usepackage{amsfonts,amssymb,stmaryrd,amscd,amsmath,latexsym,amsbsy}

\theoremstyle{plain}

\newtheorem{theorem}{Theorem}[section]

\newtheorem{thm}{Theorem}[section]

\newtheorem{prop}[theorem]{Proposition}

\newtheorem{lem}[theorem]{Lemma}

\newtheorem{cor}[theorem]{Corollary}

\newtheorem*{sol}{Solution}

\newtheorem*{conj}{Conjecture}

\theoremstyle{definition}

\newtheorem*{defi}{Definition}

\theoremstyle{remark}

\newtheorem{rema}{Remark}[section]

\newcommand{\CC}{\mathbb C}

\newcommand{\ZZ}{\mathbb Z}

\newcommand{\solu}[1]{\begin{sol}{\bf (\ref{#1})}}

\title{ Double affine Hecke algebras of rank 1 and affine cubic surfaces}

\author{ Alexei Oblomkov}

\address{Department of Mathematics, MIT, 77, Massachusetts Ave., Cambridge,
MA 02139, USA.}

\email{oblomkov@math.mit.edu}

\begin{document}

\begin{abstract}

We study the algebraic properties of the five-parameter family
$H(t_1,t_2,t_3,t_4;q)$ of double affine Hecke algebras of type
$C^\vee C_1$. This family generalizes Cherednik's double affine
Hecke algebras of rank 1.  It was introduced by Sahi and studied
by Noumi and Stokman as an algebraic structure which controls
Askey-Wilson polynomials. We show that if $q=1$, then the spectrum
of the center of $H$ is an affine cubic surface $C$, obtained from
a projective one by removing a triangle consisting of smooth
points. Moreover, any such surface is obtained as the spectrum of
the center of $H$ for some values of parameters. This result
allows one to give a simple geometric description of the action of
an extension of $PGL_2(\Bbb Z)$ by $\Bbb Z$ on the center of $H$.
When $C$ is smooth, it admits a unique algebraic symplectic
structure, and the spherical subalgebra $eHe$ of the algebra $H$
for $q=e^\hbar$ provides its deformation quantization. Using that
$H^2(C,\Bbb C)=\Bbb C^5$, we find that the Hochschild cohomology
$HH^2(H)$ (for $q=e^\hbar$) is 5-dimensional for generic parameter
values.  From this we deduce that the only deformations of $H$
come from variations of parameters. This explains from the point
of view of noncommutative geometry why one cannot add more
parameters into the theory of Askey-Wilson polynomials. We also
prove that  the five-parameter  family $H(t_1,t_2,t_3,t_4;q)$ of
algebras yields the universal deformation of the semidirect
product $q$-Weyl algebra with $\mathbb{Z}_2$ and the family of
cubic surfaces $C=C_{\underline{t}}$, $\underline{t}\in
\CC^4_{\underline{t}}$ gives the universal deformation of the
Poisson algebra $\CC[X^{\pm 1},P^{\pm 1}]^{\ZZ_2}$.

\end{abstract}

\maketitle

\section*{Introduction}

The most general family of basic hypergeometric orthogonal
polynomials are the Askey-Wilson polynomials $p_n(x;a,b,c,d;q)$
\cite{AW}. Since they were discovered, harmonic analysts tried to
find the algebraic structure which controls these polynomials.
This structure turned out to be the rank 1 double affine Hecke
algebra (DAHA) introduced by Cherednik \cite{Ch}, or, more
precisely, its generalization $H=H(t_1,t_2,t_3,t_4;q)$ to the
nonreduced root system $C^\vee C_1$, studied by Sahi, Noumi and
Stokman \cite{NS,Sa,S}. \footnote{Here the parameters
$t_1,t_2,t_3,t_4$ are related by a change of variable to
$a,b,c,d$.}

In this paper we undertake a detailed study of the algebraic
structure of $H$, and in particular discover its intimate
connection to the geometry of affine cubic surfaces.

The content of the paper is as follows. In Section 1, following
Sahi, Noumi, and Stokman, we define the double affine Hecke
algebra $H=H(t_1,t_2,t_3,t_4;q)$. If $t_i=1$, this algebra
specializes to the semidirect product \linebreak
$\Bbb C_q[X^{\pm 1},P^{\pm
1}]\rtimes\CC[\ZZ_2]$, where $\Bbb C_q[X^{\pm 1},P^{\pm 1}]$ is
the algebra of functions on the quantum torus, generated by
$X^{\pm 1},P^{\pm 1}$ with defining relation $PX=qXP$. In
particular, if $t_i=q=1$ then $H= \Bbb C[X^{\pm 1},P^{\pm
1}]\rtimes\Bbb C[\Bbb Z_2]$.

For general $\underline{t},q$, the algebra $H$ contains the
two-dimensional Hecke algebra $\Bbb C_{t_2}[\Bbb Z_2]$ of type $A_1$,
which in turn contains the symmetrizing idempotent $e$. Thus one
may define the spherical subalgebra $eHe$, which is a 5-parameter
deformation of the function algebra ${\mathcal O}(\Bbb ({\Bbb
C}^*)^2/\Bbb Z_2)$, where $\Bbb Z_2$ acts by $b\mapsto b^{-1}$.

Next we recall the Lusztig-Demazure type representation of the
algebra $H$, which is an embedding of $H(\underline{t},q)$ into an
appropriate localization of $ \Bbb C_q[X^{\pm 1},P^{\pm
1}]\rtimes\Bbb C[\Bbb Z_2]$. The existence of this representation
allows one to show that if $q=1$ then the algebra $eHe$ is
commutative and isomorphic to the center $Z(H)$ of $H$ (the Satake
isomorphism).

Finally, we recall the action of an abelian extension
$\widetilde{PGL(2,\Bbb Z)}$ of $PGL(2,\Bbb Z)$ (by $\Bbb Z$)
on $H$. More
precisely, consider the natural surjective homomorphism $\xi:
PGL(2,\Bbb Z)\to PGL(2,\Bbb Z_2)\times \Bbb Z_2=S_3\times \Bbb
Z_2$, given by $a\mapsto (a\ {\rm mod}\ 2,\det(a))$. The kernel of
$\xi$ is $\Gamma(2)/(\pm 1)$, where $\Gamma(2)$ is the level two
congruence subgroup of $SL(2,\Bbb Z)$. It turns out that the group
$\widetilde{PGL(2,\Bbb Z)}$ acts on DAHA as follows: an element
$a$ acts on the parameters of DAHA by $\xi(a)$ (where the first
component $\xi_1(a)$ acts by permuting $t_1,t_3,t_4$ and the
second component $\xi_2(a)$ by inverting the parameters).
Thus the subgroup that honestly acts on $H$ (without
changing parameters) is the preimage $\widetilde{\Gamma(2)}$ of
$\Gamma(2)/(\pm 1)$ (a subgroup of index 12).

The case $q=1$ is of special interest. In this case $q=q^{-1}$ and
the action of the abelian kernel of $\widetilde{PGL(2,\Bbb Z)}$ is
by inner automorphisms, we get an action of the group ${\rm
Ker}\xi$ on the center $Z(H)$ of $H$. Moreover, we will see that
the center of $H(\underline{t}^{-1};1)$ is canonically isomorphic
to the center of $H(\underline{t};1)$ (they have generators
satisfying the same relations), so in fact we have an action of a
twice bigger group $K={\rm Ker}\xi_1$ on the center $Z(H)$.

In Section 2, we prove that the center $Z=Z(H)$ of $H$ is
generated by three elements $X_1,X_2,X_3$ subject to one relation
$R(X_1,X_2,X_3)=0$, where $$
R(X_1,X_2,X_3)=X_1X_2X_3-X_1^2-X_2^2-X_3^2+p_1X_1+p_2X_2+p_3X_3+p_0+4,
$$ and $p_i$ are (algebraically independent) regular functions of
$t_j$. Thus the variety $C={\rm Spec}(Z)$ is an affine cubic
surface in $\Bbb C^3$, such that its projective completion $\bar
C$ differs from $C$ by three lines at infinity, forming a triangle
and consisting of smooth points. Moreover, we show that any cubic
surface with these properties is obtained in this way. This means
that a smooth projective cubic surface is obtained as $\bar C$
(generically, in 45 ways), since it contains 27 lines forming a
configuration containing 45 triangles \cite{Manin}.\footnote{We
recall that the moduli space of cubic surfaces, $PS^3\Bbb
C^4/PGL(4)$, has dimension $19-15=4$, i.e. the same as the
dimension of the moduli space of algebras $H$ with $q=1$.}.

Since for $q\ne 1$ the algebra $eHe$ is in general noncommutative,
the commutative algebra $eHe=Z(H)$ for $q=1$ has a Poisson
structure. It turns out that this Poisson structure is symplectic
at smooth points of $C$. For instance, if $C$ is smooth (generic
case) then the Poisson structure is symplectic, and is given by
the unique (up to scaling) nonvanishing 2-form $\omega$ on $C$ (one has
$H^1(C,\Bbb C)=0$ and hence any nonvanishing function on $C$ is a
constant; this guarantees the uniqueness of $\omega$). Thus the
algebra $eHe$ with $q\ne 1$ can be regarded as
a quantization of the affine cubic surface $C$ with its natural
symplectic structure.

The action of the group $K$ on $C$ has a vivid geometric
interpretation. Namely, $K$ is generated freely by three
involutions $g_1,g_2,g_3$: $$ g_1=\begin{pmatrix}1&0\\
0&-1\end{pmatrix}, g_2=\begin{pmatrix}1&2\\ 0&-1\end{pmatrix},
g_3=\begin{pmatrix}1&0\\ 2&-1\end{pmatrix}. $$ The involution
$g_i$ acts on $C$ by permuting two roots of the equation $F=0$,
regarded as a quadratic equation in $X_i$ with coefficients
depending on other variables. Thus, $g_i(X_j)=X_j$ if $j\ne i$,
and $g_i(X_i)= -X_i+X_jX_k+p_i$, $j,k\ne i$.

In Section 3, we study irreducible representations of $H$ for
$q=1$, corresponding to a given central character $\chi\in C$. We
find that the algebra $H$ is Azumaya (more specifically,
endomorphism algebra of a rank 2 vector bundle) outside of the
singular locus of $C$, while at the singular locus there are
1-dimensional representations. The singular locus is generically
empty, and always consists of at most 4 points. When there are 4
singular points (the most degenerate case), the surface is simply
the quotient of $(\Bbb C^*)^2$ by the inversion map.

In section 4 we study deformations of DAHA. We set $q=e^\hbar$,
where $\hbar$ is a formal parameter, and consider the Hochschild
cohomology of $H=H(\underline{t},q)$, in the case when
$C=C_{\underline t}$ is smooth. Since $H$ is a deformation
quantization of the function algebra $\mathcal{O}(C)$, by results
of Kontsevich the Hochschild cohomology of $H$ coincides with the
Poisson cohomology of $\mathcal{O}(C)$. On the other hand, since
$C$ is symplectic, by Brylinski's theorem, this cohomology is
equal to the De Rham cohomology of $C$. The latter is found to be:
$H^1=H^{>3}=0$, $H^2=\Bbb C^5$ (one 2-cycle on the 2-torus modulo
inversion, which is homotopically a sphere, and 4 vanishing
2-cycles at the 4 singular points of this torus). In particular,
the formal deformation space of $H$ is smooth and 5-dimensional.
The space of parameters of $H$ is also 5-dimensional. We prove
that the natural map between these spaces is generically
injective, which implies that the variation of parameters
$\underline{t},q$ produces the universal deformation of $H$. This
is to some extent an explanation why one cannot add additional
parameters into Askey-Wilson polynomials.

In section 5 we consider DAHA as deformation of the semidirect
product $\CC_q[X^{\pm 1},P^{\pm 1}]\rtimes\CC[\ZZ_2]$. We
calculate the Hochschild cohomology of the later algebra for the
case $q$ is not a root of unity. Later we show that five-parameter
family of DAHA is a universal deformation of $\CC_q[X^{\pm
1},P^{\pm 1}]\rtimes\CC[\ZZ_2]$.

Finally, in section 6 we calculate the Poisson cohomology of the
Poisson ring $\CC[X^{\pm 1},P^{\pm 1}]^{\ZZ_2}$. It turns out that
the family of cubic surfaces $\{
C_{\underline{t}}\}_{\underline{t}\in \CC^4}$ with the natural
Poisson structure gives a universal deformation of the singular
Poisson variety $Spec(\CC[X^{\pm 1},P^{\pm 1}]^{\ZZ_2})$.

 {\bf Acknowledgments.} The author thanks his adviser P. Etingof
for posing the problem, discussions, and helping to write the
introduction. The author is also grateful to J. De Jong, J. Starr
and S.L'vovski for consultations on the geometry surfaces and to
J. Stokman and V. Ostrik for their explanations about double
affine Hecke algebras. I am especially grateful to E. Rains who
helped me to correct and simplify the statements from the section
3 and find more transparent formulation of the
corollary~\ref{sings}. The author thanks referee for the useful
remarks which helped to improve the text.
 The work of the author was partially supported by the NSF grant
DMS-9988796.

\section{Double affine Hecke algebras}
\begin{defi}Let $k_0,k_1,u_0,u_1,q\in\CC^*$.
The double affine Hecke algebra $H=H(k_0,k_1,u_0,u_1;q)$ of rank
$1$ is generated by the elements $V_0$, $V_1$, $V_0^{\vee}$,
$V_1^{\vee}$ with the relations:
\begin{gather}\label{k0}
(V_0-k_0)(V_0+k_0^{-1})=0,\\
(V_1-k_1)(V_1+k_1^{-1})=0,\label{k1}\\
(V_0^{\vee}-u_0)(V^{\vee}_0+u_0^{-1})=0,\label{u0}\\
(V_1^{\vee}-u_1)(V_1^{\vee}+u_1^{-1})=0,\label{u1}\\
V_1^{\vee}V_1V_0V_0^{\vee}=q^{-1/2}.\label{VVVV}
\end{gather}
\end{defi}

This algebra is the rank one case of the algebra  discovered by
Sahi \cite{Sa}; he used this algebra to prove Macdonald's
conjectures for Koornwinder's polynomials. In our paper we use
notation from the paper \cite{NS}. In particular we use
$\underline{t}=(t_1,t_2,t_3,t_4)=(k_0,k_1,u_0,u_1)$.

\subsection{Lusztig-Demazure representation}
Let $D_q=\CC_q[X^{\pm 1},P^{\pm 1}]\rtimes\CC[\ZZ_2]$, where
$\CC_q[X^{\pm 1},P^{\pm 1}]$, $q\in \CC^*$ is a $q$-deformation of
the ring of the Laurent polynomials of two variables:
\begin{equation*}
PX=qXP.
\end{equation*}
 That is, an element of $D_q$ is a linear combination of
monomials $X^iP^js^\epsilon$, where $s\in \ZZ_2$ is the generator
of $\ZZ_2$, $i,j\in \ZZ$, $\epsilon=0,1$ and
\begin{gather*}
(X^i P^j
s)(X^{i'}P^{j'}s^{\epsilon'})=q^{-i'j}X^{i-i'}P^{j-j'}s^{\epsilon'+1}\\
(X^i
P^j)(X^{i'}P^{j'}s^{\epsilon'})=q^{i'j}X^{i+i'}P^{j+j'}s^{\epsilon'}.
\end{gather*}

In the next proposition
$D_{q,\rm loc}$ stands for the localization of $D^q$ obtained by
inverting nonzero functions of $X$.

\begin{prop}\cite{Sa,NS}\label{LDeM}
For any $q\in \CC^*$ the following formulas give an injective
homomorphism $i_q: H(k_0,k_1,u_0,u_1;q)\to D_{q,\rm loc}$:
\begin{gather*}
i_q(V_i)=T_i, i_q(V^{\vee}_i)=T^{\vee}_i,\\
T_0=k_0P^{-1}s+\frac{\bar{k}_0+\bar{u}_0X}{1-X^2}(1-P^{-1}s),\\
T_1=k_1s+\frac{\bar{k}_1+\bar{u}_1X}{1-X^2}(1-s),\\
T^{\vee}_1=X^{-1}T_1^{-1},\quad T_0^{\vee}=q^{-1/2}T_0^{-1}X.
\end{gather*}
\end{prop}

\begin{cor}\label{locLD}
If $q=1$ then $i_q$ induces an isomorphism
\begin{gather*}
H_{(\delta(X))}\simeq\CC[X^{\pm 1},P^{\pm 1}]_{(\delta(X))},\\
\delta(X)=(1-X^2)(1-k_1u_1X)(1+k_1u_1^{-1}X)
(1-k_2u_2X)(1+k_2u_2^{-1}X),
\end{gather*}
 where the element $X\in H$ is defined by $X=T_1^{-1}(T_1^\vee)^{-1}$,
and the subscript $(\delta(X))$ denotes localization by
$\delta(X)$.
\end{cor}












\subsection{Projective $PSL(2,\ZZ)$ action}
In \cite{Ch}, Cherednik showed that double affine Hecke algebras
corresponding to reduced root systems admit a projective
$SL(2,\ZZ)$-action. This action was generalized to the case of
nonreduced root systems in the papers \cite{NS,S}. Let us recall
this generalization.

We denote by $H$ the algebra
$H(\underline{t};q)=H(k_0,k_1,u_0,u_1;q)$, by $H^{\sigma}$ the
algebra $H(\sigma(\underline{t});q)=H(u_1,k_1,u_0,k_0; q)$,
$H^\tau=H(\tau(\underline{t});q)=H(u_0,k_1,k_0,u_1;q)$,
$H^\eta=H(\underline{t}^{-1};q^{-1})=
H(k_0^{-1},k_1^{-1},u_0^{-1},u_1^{-1};q^{-1})$.

\begin{prop}For any values of $q\in\CC^*$ and $\underline{t}$
\begin{enumerate}
\item
The following formulas give an isomorphism $\sigma$: $H\to
H^{\sigma}$:
\begin{gather}
\sigma(V_0)=\tilde{V}^{-1}_1\tilde{V}_1^{\vee}\tilde{V}_1,\quad
\sigma(V_1)=\tilde{V}_1,\\
\sigma(V^{\vee}_0)=\tilde{V}_0\tilde{V}_0^{\vee}\tilde{V}_0^{-1},\quad
\sigma(V^{\vee}_1)=\tilde{V}_0,
\end{gather}
where $\tilde{V}_i$, $\tilde{V}^{\vee}_i$, $i=0,1$ generate
$H^{\sigma}$ with relations (\ref{k0}-\ref{VVVV}) in which $V_i$
and $V_i^{\vee}$ are replaced by $\tilde{V}_i$ and
$\tilde{V}^{\vee}_i$ (respectively) and
$(k_0,k_1,u_0,u_1)=\underline{t}$ is replaced by
$\sigma(\underline{t})$.
\item
The following formulas give an isomorphism $\tau$: $H\to
H^{\tau}$:
\begin{gather}
\tau(V_0)=\tilde{V}_0\tilde{V}_0^{\vee}\tilde{V}_0^{-1},\quad
\tau(V_1)=\tilde{V}_1,\\ \tau(V^{\vee}_0)=\tilde{V}_0,\quad
\tau(V^{\vee}_1)=\tilde{V}_1^{\vee},
\end{gather}
where $\tilde{V}_i$, $\tilde{V}^{\vee}_i$, $i=0,1$ generate
$H^{\tau}$ with relations (\ref{k0}-\ref{VVVV}) in which $V_i$ and
$V_i^{\vee}$ are replaced by $\tilde{V}_i$ and
$\tilde{V}^{\vee}_i$ (respectively) and
$(k_0,k_1,u_0,u_1)=\underline{t}$ is replaced by
$\tau(\underline{t})$.
\item
The following formulas give an isomorphism $\eta$: $H\to H^\eta$:
\begin{gather}
\eta(V_0)=\tilde{V}_0^{-1},\quad \eta(V_1)=\tilde{V}_1^{-1},\\
\eta(V^{\vee}_0)=\tilde{V}_0(\tilde{V}_0^\vee)^{-1}\tilde{V}_0^{-1},\quad
\eta(V^{\vee}_1)=\tilde{V}_1^{-1}(\tilde{V}_1^{\vee})^{-1}\tilde{V}_1,
\end{gather}
where $\tilde{V}_i$, $\tilde{V}^{\vee}_i$, $i=0,1$ generate $H^\eta$
with relations (\ref{k0}-\ref{VVVV}) in which $V_i$ and
$V_i^{\vee}$ are replaced by $\tilde{V}_i$ and
$\tilde{V}^{\vee}_i$ (respectively) and
$(k_0,k_1,u_0,u_1)=\underline{t}$ and $q$ is replaced by
$\underline{t}^{-1}$, $q^{-1}$.
\item For all $h\in H$ we have  the relations
\begin{gather}
(\sigma\circ\sigma)(h)=V_1^{-1}hV_1,\\ (\sigma\circ\tau\circ\sigma
\circ\tau\circ\sigma\circ\tau)(h)=h,\\
(\sigma\circ\sigma\circ\tau)(h)=(\tau\circ
\sigma\circ\sigma)(h),\\ \sigma\circ\eta\circ\sigma(h)=\eta(h),
\\ \tau\circ\eta\circ\tau(h)=\eta(h).
\end{gather}
\end{enumerate}
\end{prop}

Recall that the group $PSL(2,\Bbb Z)$ has generators $$
\tilde{\sigma}=\begin{pmatrix}0&1\\-1&0\end{pmatrix},\quad
\tilde{\tau}=\begin{pmatrix}1&1\\0&1\end{pmatrix}$$ and defining
relations $$\tilde{\sigma}^2=1=(\tilde{\sigma}\tilde{\tau})^3.$$

Thus the relations from the last item of the proposition show that
the elements $\sigma,\tau$ define a projective action (i.e., an
action of a central extension) of $PSL(2,\Bbb Z)$ on the sum of
six double affine Hecke algebras (obtained from $H$ by all
permutations of $k_0,u_0,u_1$), such that the center of the
central extension preserves each of the six and acts on them by
the inner  automorphism (in particular, acts trivially on the
center of $H$).

Moreover, the group generated by the automorphisms
$\sigma,\tau,\eta$ acts on the sum of twelve double affine Hecke
algebras (obtained from $H$ by all permutations of $k_0,u_0,u_1$
and taking the inverse of all parameters including $q$). This
group is an extension (not central) of the group $PGL(2,\ZZ)$ by
$\ZZ$ acting by inner automorphisms $h\to V_1^k h V_1^{-k}$, $k\in
\ZZ$. On the center of $H$ the inner automorphisms act trivially
and we get the action of $PGL(2,\ZZ)$ on the sum of the centers of
twelve algebras. The element $\eta$ correspond to the element
$$\tilde{\eta}=\begin{pmatrix}1&0\\0&-1
\end{pmatrix}$$ of $PGL(2,\ZZ)$. In particular, here are the formulas for the involutions
$g_i$ from the introduction:
\begin{equation}\label{ggg}
g_1=\eta,\quad g_2=\eta\circ\tau^2,\quad
g_3=\eta\circ\sigma\circ\tau^2\circ\sigma.
\end{equation}

\begin{rema}
To get from the algebra $H(\underline{t};q)$ the double affine
Hecke algebra of the type $A_1$ one needs to put $k_0=u_0=u_1=1$
and $k_1=k$. The elements $\pi=V_0$, $T=V_1$,
$X=q^{1/2}V_0V_0^{\vee}$ generate this algebra with the relations:
$$ (T-k)(T+k^{-1})=0,\quad TXT=X^{-1}, \quad \pi^2=1,\quad \pi
X\pi=q^{1/2}X^{-1}.$$ These relations coincide with the relations
from the paper \cite{ChO} (see Lemma~5.7) where the double affine
Hecke algebra of the type $A_1$ was studied carefully. Let us
remark that under this degeneration the isomorphisms $\sigma$,
$\tau$ become automorphisms of this algebra, so we get a
projective action of $PSL(2,\Bbb Z)$ on $H$.
\end{rema}







\section{Affine cubic as the spectrum of the
center of $H(\underline{t};1)$}
Now we restrict ourselves to the case $q=1$. We denote
$H(\underline{t};1)$ by $H$. In this section we prove the
following theorem, which is one of our main results.

\begin{thm}\label{eqnCub}
For any values $k_0,k_1,u_0,u_1$ the elements
\begin{gather}\label{frml cntr}
X_1=V_1^{\vee}V_1+V_0V_0^{\vee}, \quad
X_2=V_1V_0+V_0^{\vee}V_1^{\vee},\quad
X_3=V_1V_0^{\vee}+(V_0^{\vee})^{-1}V_1^{-1},
\end{gather}
 generate the center
$Z(H)$ of the double affine Hecke algebra $H=H(\underline{t};1)$.
Moreover we have $$Z(H)=\CC[X_1,X_2,X_3]/(R_{\underline{t}}),$$
\begin{multline*}
R_{\underline{t}}=X_1X_2X_3-X_1^2-X_2^2-X_3^2+
(\bar{u}_0\bar{k}_0+\bar{k}_1\bar{u}_1)X_1+
(\bar{u}_1\bar{u}_0+\bar{k}_0\bar{k}_1)X_2+\\ (\bar{k}_0\bar{u}_1+
\bar{k}_1\bar{u}_0)X_3+\bar{k}_0^2+\bar{k}_1^2+\bar{u}_0^2+\bar{u}_1^2-
\bar{k}_0\bar{k}_1\bar{u}_0\bar{u}_1+4,
\end{multline*}
where $\bar{k}_i=k_i-k_i^{-1}$, $\bar{u}_i=u_i-u_i^{-1}$.
\end{thm}

The proposition is proved in subsection \ref{pfpropo}.

\subsection{Properties of the affine cubic surface}
Before proving this proposition let us list the properties of the
affine cubic $C_{\underline{t}}=\{R_{\underline{t}}(X)=0\}$ and
its completion $\bar{C}_{\underline{t}}\subset\mathbb P^3$. We use
the term ``triangle'' for the union  of three distinct lines in
the plane, which don't intersect in the same point.
\begin{prop}
\begin{enumerate}
\item For any values of the parameters $\underline{t}$ the complement
$\bar{C}_{\underline{t}}\setminus C_{\underline{t}}$ is a
triangle.
\item For any value of $\underline{t}$ the cubic
$\bar{C}_{\underline{t}}$ is irreducible, normal, has a finite
number of the  singular points, and is smooth at infinity.
\item If $C'$ is an affine irreducible cubic with a triangle at
infinity consisting of smooth points, then $C'=C_{\underline{t}}$
for some $\underline{t}\in\CC^4$.
\item If $\underline{t}$ is generic then $C_{\underline{t}}$ is
smooth.
\end{enumerate}
\end{prop}
\begin{proof}
{\it (1)} If we rewrite the equation for $C_{\underline{t}}$ in terms
of the homogeneous coordinates $X_i=x_i/x_0$ then it is easy to
see that the intersection of the infinite plane with
$C_{\underline{t}}$ is given by the equations:
$$
x_1x_2x_3=0,x_0=0.
$$

{\it (2)} The irreducibility of $\bar C_{\underline{t}}$ and
smoothness at infinity are immediate from the equation. The
finiteness of the number of singular points follows from
smoothness at infinity. To prove normality, note that $\bar
C_{\underline{t}}$ is a Cohen-Macaulay variety. Since it has a
finite number of the singular points, by Serre's criterion it is
normal.

{\it (3)} Suppose  that the coordinates are chosen in such a way
that the three lines lying on the infinite plane $\{ x_0=0\}$ are
given by the equations $x_0=x_1x_2x_3=0$. Then the equation of the
cubic surface has the form
$$X_1X_2X_3=\sum_{i,j}a_{ij}X_iX_j+\sum_k b_k X_k+c.$$ Making the
shifts $X_i\to X_i+2a_{jk}$, $j,k\ne i$ and rescaling $X_i$ we
arrive at the equation: $$X_1X_2X_3=\sum_{i=1}^3\epsilon_i
X_i^2+\sum_k b'_kX_k+c',$$ where $\epsilon_i=0,1$. If
$\epsilon_i=0$ then the infinite point with homogeneous
coordinates $x_i=1,x_j=0, j\ne i$ is a singular point of $C'$. So
the proof follows from Lemma~\ref{pi} which we prove in
section~\ref{frep}.

{\it (4)} The last item follows from the previous one because the
generic cubic surface is smooth and contains a triangle of lines.
\end{proof}

\subsection{The action of $PGL(2,\ZZ)$ on the cubic surface}
We can easily calculate the action of $PGL(2,\ZZ)$ on the elements
$X_i$:
\begin{gather}\label{SL2centerF}
\sigma(X_1)=X_2^\sigma,\quad\sigma(X_2)=X_1^\sigma\\
\sigma(X_3)=X_1^\sigma
X_2^\sigma-X_3^\sigma+\bar{u}_0\bar{u}_1+\bar{k}_0\bar{k}_1\\
\tau(X_1)=X_1^\tau,\quad\tau(X_3)=X_2^\tau\\ \tau(X_2)= X_1^\tau
X_2^\tau-X^\tau_3+\bar{u}_0\bar{u}_1+\bar{k}_0\bar{k}_1,\\
\eta(X_1)=X_1^\eta,\quad \eta(X_2)=X_2^\eta\\ \eta(X_3)=X^\eta_1
X^\eta_2-X^\eta_3+\bar{u}_0\bar{k}_1+\bar{k}_0\bar{u}_1\label{SL2centerL}
\end{gather}
where $X_i^\sigma$, $X_i^\tau$, $X^\eta_i$ are the corresponding
generators of the center of $H^\sigma$, $H^\tau$, $H^\eta$
respectively.
Using formulas
(\ref{ggg}) we get the formulas for the action of  $g_i$ from the
introduction.

\begin{rema}\label{foreqnCub} The formulas
(\ref{SL2centerF}-\ref{SL2centerL}) define the regular maps
$\sigma,\tau,\eta$: $\mathbb {A}^3\to \mathbb{A}^3$. The direct
calculation shows that $\sigma,\tau,\eta$ map the surface
$C_{\underline{t}}\subset\mathbb{A}^3$ into the surfaces
$C_{\sigma(\underline{t})}$, $C_{\tau(\underline{t})}$,
$C_{\eta(\underline{t})}$, respectively.
\end{rema}

\begin{cor} Formulas (\ref{SL2centerF}-\ref{SL2centerL}) yield an
embedding of the congruence subgroup  $K$ into  the group
$Aut(C_{\underline{t}})$ of automorphisms of the affine cubic
$C_{\underline{t}}$.
\end{cor}
\begin{proof}
The existence of the map follows from the previous reasoning. The
injectivity of the map is an open condition on the parameters
$\underline{t}$. Hence it is enough to prove that it is an
injection for some particular value of $\underline{t}$.

If $\bar{k}_0=\bar{k}_1=\bar{u}_0=\bar{u}_1=0$ then algebra
$H(\underline{t};1)$ becomes $\CC[X^{\pm 1},Y^{\pm 1}]\rtimes
\ZZ_2$ and $\CC[X^{\pm 1},Y^{\pm 1}]^{\ZZ_2}$. The action of $K$
on $\CC[X^{\pm 1},Y^{\pm 1}]^{\ZZ_2}$ is induced by the standard
action of $GL(2,\ZZ)$ on $\CC[X^{\pm 1},Y^{\pm 1}]$: $X^iY^j\to
X^{\gamma(i)}Y^{\gamma(j)}$, $\gamma\in GL(2,\ZZ)$. This action
obviously induces an injection: $K\to Aut(C_{\underline{t}})$.
\end{proof}

\begin{rema}
The automorphisms from the corollary do not extend to
automorphisms of the projective cubic surface. Their extensions
are birational automorphisms of the projective surface.
 Moreover,
the formulas (\ref{SL2centerF}-\ref{SL2centerL}) extend the action
of $K$ to the reducible projective surface $C_{\underline{t}}\cup
L$, where $L$ is infinite plane. The group $K$ acts on the
infinite plane $L$ by Cremona transformations.
\end{rema}

\subsection{Proof of Theorem~\ref{eqnCub}}\label{pfpropo}
From  Proposition~\ref{LDeM} it follows that the element $X_1$ is
central. The elements $X_2$, $X_3$ are central because they are
the results of the $PGL(2,\ZZ)$ action on $X_1$.

In the Appendix we prove that $X_i$ satisfy the cubic equation
$R_{\underline{t}}(X)=0$. It is a routine but rather long
calculation.

We now prove that $X_1,X_2,X_3$ generate the center $Z$. As the center
of $\CC[X^{\pm 1},P^{\pm 1}]\rtimes \ZZ_2$ is equal to $\CC[X^{\pm
1},P^{\pm 1}]^{\ZZ_2}$, the corollary~\ref{locLD} implies
that the map $i_1$ induces a birational isomorphism between
$Spec(Z)$ and $Spec(\CC[X^{\pm 1},P^{\pm 1}]^{\ZZ_2})$. The ring
$\CC[X^{\pm 1},P^{\pm 1}]^{\ZZ_2}$ is generated by the elements $$
I_1=X+X^{-1},\quad I_2=P+P^{-1},\quad I_3=XP+P^{-1}X^{-1}$$ modulo
the relation $$ I_1I_2I_3=I_1^2+I_2^2+I_3^2-4.$$ An easy
calculation shows that
\begin{gather*}
X_1=I_1,\\ X_2(I_1^2-4)=Q_0(I_1)+Q_2(I_1)I_2+Q_3(I_1)I_3,\\
X_3(I_1^2-4)=S_0(I_1)+S_2(I_1)I_2+S_3(I_1)I_3,
\end{gather*}
where $\deg Q_2=2$ with the leading coefficient $k_1/k_0$, $\deg
S_3=2$ with the leading coefficient $k_0/k_1$, $\deg Q_3\le 1$,
$\deg S_2\le 2$.

Hence we can express $I_2,I_3$ through $X_i$:
\begin{gather*}
I_2=(X_2R_2^2(X_1)+X_3R_2^3(X_1)+R_2^0(X_1))/D_1(X_1),\\
I_3=(X_2R_3^3(X_1)+X_3R_3^3(X_1)+R_3^0(X_1))/D_1(X_1),
\end{gather*}
where $R^i_j$ are polynomials and $D_1=Q_2S_3-Q_3S_2$ is a
polynomial of degree $4$ with the leading coefficient $1$. As a
consequence we get a birational isomorphism between $Spec(Z)$ and
the cubic $C_{\underline{t}}$, under which the image of any
element $F$ of the center $Z$ can be presented in the form
$F=R_1(X)/(D_1(X_1)(X_1^2-4))^N$, where $R_1$ is a polynomial.

As we have a (projective) $PSL(2,\ZZ)$ action on the DAHA,
recalling remark~\ref{foreqnCub} we see that the same statement
holds for $X_2$ and $X_3$. That is, any element $F\in Z$ has a
representation in the form $F=R_i(X)/(D_i(X_i)(X_i^2-4))^{N_i}$
where $i=1,2,3$ and  $R_i$ are polynomials depending on $F$, while
$D_i$ is a fixed polynomial. In other words any regular function
on $Spec(Z)$ is a rational function on $C_{\underline{t}}$ with
singularities on a set of codimension $2$. But we know that
$C_{\underline{t}}$ is normal, hence any function with
singularities at codimension $2$ is regular. Thus we proved that
$Z= \CC[X_1,X_2,X_3]/(R_{\underline{t}})$.


\section{Finite dimensional representations of $H$ and the
spherical subalgebra}\label{frep} Let $e=(1+k_1V_1)/(1+k_1^2)$ be
the symmetrizer. In all places where we use the symmetrizer  we
 suppose that $k_1^2\ne -1$ and we denote by $C_{k_1}[\mathbb{Z}_2]$
 the subalgebra of $H$ generated by $V_1$. The
algebra $eHe$ is called the {\it spherical subalgebra} of $H$. The
space $He$ has the structure of a right $eHe$-module.

Theorem 5.1 and Corollary 6.2   from \cite{O} say
\begin{prop}\label{fromO}
\begin{enumerate}
\item The map $z\mapsto z e$ is an isomorphism $Z\to eHe$. \item
The left action of $H$ on $He$ induces  an isomorphism of algebras
$H\simeq End_{eHe}(He)$. \item If $C_{\underline{t}}$ is smooth
then $He$ is a projective $eHe$-module which corresponds to a
vector bundle of rank $2$ on $C_{\underline{t}}$. \item If
$U\subset C_{\underline{t}}$ is a smooth affine open set and $E$
is an irreducible representation of $H$ with central character
$\chi$ belonging to $U$ then $E$ is the regular representation of
$\mathbb{C}_{k_1}[\ZZ_2]$, namely $$ E=He\otimes_{eHe}\chi.$$
\end{enumerate}
\end{prop}

\begin{rema}
In \cite{O}, we treat the case of reduced root systems, but the
proof is the same in the non-reduced case.
\end{rema}

It turns out that the statement converse to the last part of the
proposition~\ref{fromO} holds, and we can give a simple
representation-theoretic description of the locus of the singular
affine cubic surfaces. For that we need to understand better the
structure of the map $p$: $\mathbb C_{\underline{t}}^4\to\mathbb
C^4_p$, given by the coefficients of the cubic surface
$C_{\underline{t}}$:
\begin{gather*}
p_1=\bar{u}_0\bar{k}_0+\bar{k}_1\bar{u}_1,\quad
p_2=\bar{u}_1\bar{u}_0+\bar{k}_0\bar{k}_1,\\
p_3=\bar{k}_0\bar{u}_1+ \bar{k}_1\bar{u}_0,\quad p_0=
\bar{k}_0^2+\bar{k}_1^2+\bar{u}_0^2+\bar{u}_1^2-
\bar{k}_0\bar{k}_1\bar{u}_0\bar{u}_1.
\end{gather*}  The structure of the map becomes transparent if one
introduce  torus $\mathbb{T}$ and maps $\theta$: $\mathbb
C^4_{\underline{t}}\to\mathbb{T}$, $\pi$: $\mathbb{T}\to \mathbb
C^4_p$ with the property $p=\pi\circ\theta$. Now let us introduce
these maps and torus $\mathbb{T}$.

Let $\mathbb{T}=Spec(\mathbb C[P])$ where $P$ is the weight
lattice for the root system $D_4$. By the definition
$$P=\oplus^4_{i=1} \mathbb Z\varepsilon_i+\mathbb Z(\frac12
(\varepsilon_1+\varepsilon_2+\varepsilon_3+\varepsilon_4)),$$ and
we can think about this algebra as the algebra generated by
$s_1^{\pm 1}$, $s_2^{\pm 1}$, $s_3^{\pm 1}$, $s_4^{\pm 1}$,
$\delta^{\pm 1}$ modulo the relation $s_1s_2s_3s_4=\delta^2$.
There is a natural action of the Weyl group $W=W_{D_4}$ on
$\mathbb{T}$.

We can embed the algebra $\mathbb C[P]$ into the algebra $\mathbb
C[t_1^{\pm 1},t_2^{\pm 1},t_3^{\pm 1},t_4^{\pm 1}]$ by the
formulas:
$$ s_1=t_1t_2,\quad s_2=-t_1/t_2,\quad s_3=-t_3/t_4,\quad
s_4=t_3t_4,\quad \delta=t_1t_3. $$ Let us denote the corresponding
map $\mathbb C^4_{\underline{t}}\to \mathbb{T}$ by $\theta$.

The direct calculation shows that we have the decomposition
$p=\pi\circ \theta$ for the map $p$: $\mathbb
C^4_{\underline{t}}\to \mathbb C^4_{p}$ with $\pi$: $\mathbb{T}\to
\mathbb C^4_p$ given by the formulas:
$$p_0=m_{\omega_2}(s)-4,\quad p_1=m_{\omega_4}(s),\quad
p_2=m_{\omega_1}(s),\quad p_3=m_{\omega_3}(s),$$ where $\omega_i$
are the fundamental weights:
$$ \omega_1=\varepsilon_1,\quad
\omega_2=\varepsilon_1+\varepsilon_2,\quad
\omega_3=\varepsilon_1+\varepsilon_2+\varepsilon_3-\varepsilon_4,
\quad
\omega_4=\varepsilon_1+\varepsilon_2+\varepsilon_3+\varepsilon_4,$$
and $m_{\omega_i}(s)=\sum_{\lambda\in W\omega_i}s^\lambda$ are the
orbit sums. As the orbit sums $m_{\omega_i}$, $i=1,2,3,4$ freely
generate the ring $\mathbb{C}[P]^W$ of $W$-invariants we get

\begin{lem}\label{pi}
The map $\pi$: $\mathbb{T}\to \mathbb{C}_p$ is an epimorphism and
for any $s\in\mathbb{T}$ we have $\pi^{-1}(\pi(s))=\cup_{w\in W}
w(s)$.
\end{lem}

\begin{rema} The map $\theta$ is the Galois covering with the
Galois group $\mathbb{Z}_2$. \end{rema}

Let us denote by $C_{\underline{s}}$ the surface
$C_{\underline{t}}$ with
$\underline{t}\in\theta^{-1}(\underline{s})$

\begin{prop}
 The point $\chi\in C_{\underline{s}}$ is singular if and only if
there exists an element $w\in W$ such that the one dimensional
representation $\chi$ of $Z=Z(H(\underline{t},1))$,
$\underline{t}\in\theta^{-1}(w(s))$ extends to a one dimensional
representation of $H(\underline{t};1)$.
\end{prop}
\begin{proof}
Proposition~\ref{fromO} implies that if $\rho$:
\begin{equation*}
V_0\mapsto t_1,\quad V_1\mapsto t_2,\quad V_0^{\vee}\mapsto
t_3,\quad V_1^{\vee}\mapsto t_4,
\end{equation*}
is an extension of the one-dimensional representation $\chi$:
\begin{equation*}
\chi=(t_2t_4+t_1t_3,t_1t_2+t_3t_4,t_2t_3+t_2^{-1}t_3^{-1}),
\end{equation*}
of $Z$ then the point $\chi\in C_{\underline{t}}$ is a singular
point. So we need to prove the converse.

It is well known that variety $S\subset\CC^4_p$ of the singular
cubic surfaces is irreducible. By the first part of the
proposition the map $\pi$: $\underline{s}\mapsto
\pi(\underline{s})$ sends $\Sigma=\{ s_1s_4=1\}$ onto $S$.
Obviously $\pi^{-1}(S)\subset\CC_{p}^4$ is a union of the
hypersurfaces. From the Lemma~\ref{pi} we get that $\pi^{-1}(S)=
\cup_{w\in W} w(\Sigma)$.
\end{proof}

Let us denote by $\Sigma\subset\CC^4_{\underline{s}}$ the locus of
the points $\underline{s}$ with the property that the surface
$C_{\underline{s}}$ has at least one singular point; by
$\Sigma'\subset\Sigma$ the locus of the surfaces with at least two
singular points (counted with multiplicities); by
$\Sigma''\subset\Sigma'$ the locus of the surfaces with at least
three singular points (counted with multiplicities); by
$\Sigma'''\subset\Sigma''$ the locus of the surfaces with at least
four singular points (counted with multiplicities). By the
multiplicity of a singular point we mean the  Milnor number
\cite{Arn}. Having in the mind this definition we get the
following result.

\begin{cor}\label{sings}
The surface $C_{\underline{s}}$ can have only $ADE$ singularities.
The type of the singularity of the surface $C_{\underline{s}}$ is
the type of the stabilizer $Stab(s)\subset W$. That is the list
below gives the complete classification of the possible
singularities of the surface $C_{\underline{s}}$:
\begin{enumerate}
\item $\Sigma=\cup_{1\le i<j\le 4,\epsilon=\pm
1}\Sigma^\epsilon_{ij}$ and $\underline{t}\in\Sigma_{ij}^\epsilon$
if and only if $s_i^\epsilon=s_j$. If
$\underline{s}\in\Sigma\setminus\Sigma'$ then $C_{\underline{t}}$
has one singular point of type $A_1$.

\item $\Sigma'=\Sigma'_2\cup\Sigma'_{1,1}$, where
$$\Sigma'_2=\cup_{\{i,j,k,l\}=\{1,2,3,4\},\epsilon=\pm
1}\{s_i=s_j=s_k^\epsilon\},$$
$$\Sigma''_{1,1}=\hat{\Sigma}''_{1,1}\cup\tilde{\Sigma}''_2,$$
$$\hat{\Sigma}''_{1,1}=\cup_{i<j,\epsilon=\pm 1}
\{ s_i=\epsilon,s_j=\epsilon\},$$
$$\tilde{\Sigma}''_{1,1}=\cup_{\{i,j,k,l\}=\{1,2,3,4\},\epsilon\in
\{\pm 1\}^2}\{s_i=s_j^{\epsilon_1},s_k=s_l^{\epsilon_2}\},$$ and
if $\underline{s}\in\Sigma'_2\setminus\Sigma''$ then the singular
locus of $C_{\underline{s}}$ consists of one point of type $A_2$;
if $\underline{s}\in\Sigma'_{1,1}\setminus \Sigma''$ then the
singular locus consists of two distinct points of type $A_1$.

\item $\Sigma''=\Sigma''_3\cup\Sigma''_{1,1,1}$, where
$$\Sigma''_3=\hat{\Sigma}''_3\cup\tilde{\Sigma}''_{1,1,1},$$
$$
\hat{\Sigma}''_3=\cup_{\epsilon\in\{\pm 1\}^2}\{s_1^{\epsilon_1}=
s_2^{\epsilon_2}=s_3^{\epsilon_1\epsilon_2}= s_4\},$$
$$\tilde{\Sigma}''_3=\cup_{\{i,j,k,l\}=\{1,2,3,4\},\epsilon=\pm 1}\{
s_i=s_j=s_k=\epsilon\},$$
\begin{equation*}\Sigma''_{1,1,1}=\cup_{\{ i, j, k,
l\}=\{1,2,3,4\}, \epsilon\in\{\pm 1\}^2}\{s_i=s_j=\epsilon_1,
s_k^{\epsilon_2}=s_l\},
\end{equation*}
 and if
$\underline{s}\in \Sigma''_3\setminus\Sigma'''$ then the singular
locus of $C_{\underline{s}}$ consists of one point of type $A_3$;
if $\underline{s}\in\Sigma''_{1,1,1}\setminus\Sigma'''$ then the
singular locus consists of three distinct points of type $A_1$.

\item $\Sigma'''=\Sigma'''_4\cup\Sigma'''_{1,1,1,1}$ where
$$\Sigma'''_4=\cup_{\epsilon\in\{\pm
1\}} \{s_1=s_2=s_3=s_4=\epsilon\}.$$
$$\Sigma'''_{1,1,1,1}=
\cup_{\{i,j,k,l\}=\{1,2,3,4\}}\{s_i=s_j=-s_k=-s_l\},$$ and if
$\underline{s}\in\Sigma'''_4$ then the singular locus of
$C_{\underline{s}}$ consists of one point of type $D_4$; if
$\underline{s}\in\Sigma'''_{1,1,1,1}$ then the singular locus of
$C_{\underline{s}}$ consists of four distinct points of type
$A_1$.
\end{enumerate}
\end{cor}
\begin{proof}

 The fact that isolated singularities of the cubic surface which
is not a cone are  $ADE$ is very classical \cite{Caley}.

Items (1)-(3) are simple corollaries of the  classification of
$ADE$ diagrams, because in all cases it is easy to calculate the
Milnor number of the singular point and the Milnor number is equal
to the number of nodes in the diagram.

Let us give  more details for the items (1), (2). The rest items
are absolutely analogous. First of all Proposition~\ref{sings}
implies that  if $\underline{s}\in\Sigma^{1}_{ij}$ then the point
$\chi\in C_{\underline{s}}$,
$\chi=(\delta(1+1/(s_ks_l)),-s_i-s_j,\delta/s_i(1+1/(s_ks_l)))$
with $\{i,j,k,l\}={1,2,3,4}$ is
 a singular point. Moreover for the point of $\Sigma^1_{ij}$ which does
 not belong to any $\Sigma^{\epsilon}_{i'j'}$, $(\epsilon,i',j')\ne
 (1,i,j)$ the point $\chi$ is the only singular
 point of $\Sigma_{s}$ and it has
type $A_1$. Analogously, for the generic (in the same sense as
before) $\underline{s}\in\Sigma^{-1}_{ij}$ the point
$\chi=(\delta(1+1/(s_ks_l)),-s_i-s_j,\delta(1/s_l+1/s_k))$ is the
only singular point of $C_{\underline{s}}$ and it is of  type
$A_1$.

To prove (2) we need to study singularities of the surfaces
$C_{\underline{s}}$ with
$\underline{s}\in\Sigma^{\epsilon}_{ij}\cap\Sigma^{\epsilon'}_{i'j'}$
, $(\epsilon,i,j)\ne(\epsilon',i',j')$ such that $\underline{s}$
does not belong to any triple intersection of the surfaces
$\Sigma^{\epsilon}_{ij}$. Thus  the formulas from the previous
paragraph imply that
 $C_{\underline{t}}$ has one singular point if
$\{i,j\}\cap\{k,l\}\ne \emptyset$ and two singular points if
$\{i,j\}\cap\{k,l\}=\emptyset$ or $\{i,j\}=\{k,l\}$. The
topological interpretation of the Milnor number \cite{Arn} implies
that in the first case the Milnor number of the singularity is
equal to $2$.

In item (4), we have two candidates for the type of singularities
of $C_{\underline{s}}$, $\underline{s}\in \Sigma_4$: $A_4$ and
$D_4$. We can distinguish them by calculating the rank of the
Hessian which is an invariant of the singularity. In our case the
rank of the Hessian is equal to $1$, that is, we have the
singularity of type $D_4$.
\end{proof}

\begin{cor} The surface $C_{\underline{t}}$ has no more then
four  singular
points and it has four singular points if and only if
$H(\underline{t};1)=\CC[X^{\pm 1},P^{\pm 1}]\rtimes\CC[\ZZ_2]$
\end{cor}

\begin{rema} The four-parameter family of  the cubic surfaces
$\{C_p\}_{p\in\CC_p^5}$ gives the miniversal deformation of the
singularity $D_4$ by the criterion from the page 51 of the book
\cite{Arn}.
\end{rema}

\section{Hochschild cohomology of $H(\underline{t};q)$ and
noncommutative deformations of the affine cubic}

\subsection{Topology of the affine cubic surface}
 Let us calculate the cohomology of $C_{\underline{t}}$:
\begin{prop} If $C_{\underline{t}}$ is smooth
then
\begin{enumerate}
\item
$H^i(C_{\underline{t}},\CC)=0$ if $i=1$, or $i>2$
$H^0(C_{\underline{t}},\CC)=\CC$,
$H^2(C_{\underline{t}},\CC)=\CC^5$.
\item Any algebraic regular nonvanishing function on
$C_{\underline{t}}$ is a constant.
\end{enumerate}
\end{prop}
\begin{proof}
Applying the weak Lefschetz theorem
to $\overline{C}_{\underline{t}}$ and the infinite hyperplane in $\Bbb P^3$,
we find that $H^{>2}(C_{\underline{t}},\CC)=H^1(C_{\underline{t}},\CC)=0$.
On the other hand, it is well known that
$\overline{C}_{\underline{t}}$ is the blow
up of the projective plane in six points. Hence
$\chi(\overline{C_{\underline{t}}})=9$. On the other hand we know
that $\overline{C_{\underline{t}}}\setminus C_{\underline{t}}$ is
a triangle, and its Euler characteristic is equal to $3$. Hence
$\chi(C_{\underline{t}})=6$ and we proved the first part of the
statement.

The second part of the statement follows from the first one.
Indeed, if we have a nonconstant nonvanishing function $f$ on
$C_{\underline{t}}$ then $d\ln(f)$ is a regular differential. As
$H^1(C_{\underline{t}},\CC)=0$, by Grothendieck's theorem $\ln(f)$
an algebraic function. It is impossible since it has logarithmic
growth.
\end{proof}

\subsection{Symplectic structure on the affine cubic surface:
connection with noncommutative deformation of $H$}
Assume again that $C_{\underline{t}}$ is smooth.
 As the ratio of two regular nonvanishing $2$-forms is a
nonvanishing function, we have
\begin{cor} Any algebraic symplectic form on $C_{\underline{t}}$
is proportional to the form: \begin{equation} \label{2form}
\Omega=dX_1\wedge dX_2/\frac{\partial R_{\underline{t}}}{\partial
X_3}.\end{equation}
\end{cor}

The symplectic form $\Omega$ yields a nondegenerate Poisson
bracket on $\CC[C_{\underline{t}}]$:
\begin{gather}
\{X_1,X_2\}=X_1X_2-2X_3+\bar{k}_0\bar{u}_1+\bar{k}_1\bar{u}_0,\\
\{X_2,X_3\}=X_2X_3-2X_1+\bar{k}_0\bar{u}_0+\bar{k}_0\bar{u}_0,\\
\{X_3,X_1\}=X_3X_1-2X_2+\bar{u}_0\bar{u}_1+\bar{u}_0\bar{u}_1.
\end{gather}

On the other hand  $eH(q,\underline{t})e$, $q=e^{h}$  gives a
natural noncommutative deformation of $\CC[C_{\underline{t}}]$
which induces the Poisson bracket
$$\{F,G\}_{def}=[F_h,G_h]/h|_{h=0},$$ where $F_h,G_h\in
H\otimes\CC[[h]]$ have the property $F_0=F$, $G_0=G$. Moreover
these two brackets coincide:

\begin{prop}\label{Poisbr} For any values $\underline{t}$ we have
$$ \{\cdot,\cdot\}=\{\cdot,\cdot\}_{def}.$$
\end{prop}
\begin{proof}
It suffices to prove the result for those $\underline{t}$ for which
$C_{\underline{t}}$ is smooth.

If we show that the form $\Omega_{def}$ corresponding to the
Poisson bracket is a regular algebraic nonvanishing form, then we show that
$\{\cdot,\cdot\}$ is proportional to $\{\cdot,\cdot\}_{def}$.
Let us do so.

 Obviously, the Poisson bracket  $\{\cdot,\cdot\}_{def}$ is regular,
hence $\Omega_{def}$ has no zeroes. To see that
$\{\cdot,\cdot\}_{def}$ is nondegenerate it is enough to prove
this for the set of codimension $2$ at $C_{\underline{t}}$ because
$C_{\underline{t}}$ is smooth. By Proposition~\ref{locLD} the
bracket $\{\cdot,\cdot\}_{def}$ coincides on the open subset
$U=\{X,\delta(X)\ne 0\}$ with the standard Poisson bracket on
$\CC[X^{\pm 1},P^{\pm 1}]^{\ZZ_2}$ which is nondegenerate. The
same is true for the open subset $\gamma(U)$, $\gamma\in K$. But
the complement to the covering by the sets $\gamma(U)$ has
codimension $2$.

The coefficient of proportionality can be easily calculated on the
described open subset $U$.

In the Appendix we give a direct proof of this proposition.
\end{proof}

\subsection{Hochschild cohomology of $H(\underline{t};q)$}
 We denote by $\mathcal{H}_{\underline{t}}$
the algebra $H(\underline{t};q)$, $q=e^h$, considered as
topological $\CC[[h]]$ module (that is, the parameter $h$ is formal).

In this subsection we consider only $H(\underline{t};q)$ with
$C_{\underline{t}}$ smooth.

Let $\mathcal{K}=\CC((h))$ and
$\mathcal{H}_{\mathcal{K}}=\mathcal{H}_{\underline{t}}
\hat{\otimes}_{\CC[[h]]}\mathcal{K}$ (a $\mathcal{K}$-algebra).
Here $\hat{\otimes}$ stands for the completed tensor product.

We will need the following corollary of Kontsevich's quantization
theory \cite{Ko1,Ko2}.

 Let $A$ be a commutative associative $\CC$-algebra, and
$\mathcal{A}$ a topologically free $\CC[[h]]$-algebra such that
$\mathcal{A}/h\mathcal{A}\simeq A$. Let
$\mathcal{A}_{\mathcal{K}}=\mathcal{A}\hat{\otimes}_{\CC[[h]]}\mathcal{K}$
be the completed tensor product. The deformation $\mathcal{A}$ of
$A$ gives rise to a Poisson structure $\{\cdot,\cdot\}_{def}$ on
$A$, making $Spec(A)$ a Poisson variety.

\begin{prop}\label{fromEG} Assume that the canonical isomorphism:
$\mathcal{A}/h\mathcal{A}\simeq A$ can be lifted to a topological
$\CC[[h]]$-module isomorphism: $\mathcal{A}\simeq A[[h]]$. If
$Spec(A)$ is smooth and the Poisson structure
$\{\cdot,\cdot\}_{def}$ on $A$ is non-degenerate, i.e. makes
$Spec(A)$ a symplectic manifold, then there is a graded
$\mathcal{K}$-vector space isomorphism:
$HH^{\bullet}(\mathcal{A}_{\mathcal{K}})\simeq\mathcal{K}\otimes_{\CC}
H^{\bullet}(Spec(A),\CC)$.
\end{prop}

\begin{rema} This statement can be generalized to the case when
$Spec(A)$ is a smooth Poisson variety. In this case the cohomology
of $Spec(A)$ should be replaced with the Poisson cohomology of
$A$. This generalization of the proposition follows from Formality
theory (see Proposition 15.2 in \cite{EG}).
\end{rema}

The last item of Proposition~\ref{fromO} implies that algebra
$e\mathcal{H}_{\underline{t}}e$ is Morita equivalent to the
algebra $\mathcal{H}_{\underline{t}}$. Hence application of
proposition~\ref{fromEG} yields
\begin{cor} If the cubic surface $C_{\underline{t}}$ is smooth then
$$
HH^{\bullet}(e\mathcal{H}_{\mathcal{K}}e)=
HH^{\bullet}(\mathcal{H}_{\mathcal{K}})\simeq
H^{\bullet}(C_{\underline{t}},\CC)\otimes_{\CC}\mathcal{K}.$$
\end{cor}

We are not able to calculate $HH^{\bullet}(H(\underline{t};q))$
but the previous statement gives a clue for the possible answer.
\begin{conj} For the generic values $q,\underline{t}$ we have:
$$HH^1(H(\underline{t};q))=HH^{>2}(H(\underline{t};q))=0,\quad
HH^2(H(\underline{t};q))=\CC^5.$$
\end{conj}

\subsection{Deformations of the DAHA}
In this subsection we prove
\begin{thm}\label{maindef}For a generic value of $\underline{t}$, the family
$\{e\mathcal{H}_{\underline{s}} e\}_{\underline{s}\in \CC^4}$
gives a universal deformation of algebra
$e\mathcal{H}_{\underline{t}}e$, and the family
$\{\mathcal{H}_{\underline{s}}\}_{\underline{s}\in\CC^4}$ gives a
universal deformation of $\mathcal{H}_{\underline{t}}$.
\end{thm}
\begin{proof}
The second Hochschild cohomology
$HH^2(\mathcal{H}_{\mathcal{K}})$ describes the first order
deformations of $\mathcal{H}_{\mathcal{K}}$. As
$HH^3(\mathcal{H}_{\underline{t}})=0$, all first order
deformations are unobstructed, and the tangent space to the space
of deformations of $\mathcal{H}_{\underline{t}}$ is equal to
$\mathcal{K}^5$.

On the other hand to any element
$(f(h),\underline{s}(h))\in(\CC[[h]])^5$ corresponds the formal
deformation $H(e^{h(1+\varepsilon f(h))},\underline{t}+
\varepsilon\underline{s}(h))$ of
$\mathcal{H}_{\underline{t}}$. That is, we have a morphism of
$\CC[[h]]$-modules $\Phi_{\underline{t}}$: $(\CC[[h]])^5\to
HH^2(\mathcal{H}_{\underline{t}})/({\rm torsion})$,
and to prove the theorem we
need to prove that this map is injective. Obviously
$\Phi_{\underline{t}}$ is injective if the induced
map $\Phi'_{\underline{t}}:(\CC[[h]])^5/h(\CC[[h]])^5\simeq\CC^5\to
H^2(C_{\underline{t}})\simeq HH^2(\mathcal{H}_{\underline{t}})/ [h
HH^2(\mathcal{H}_{\underline{t}})+({\rm torsion})]$ is injective.

Let us observe that in the space $H^2(C_{\underline{t}},\CC)$
there is a special element $\Omega$ corresponding to the
symplectic form.

 From the construction of the isomorphism
$HH^2(\mathcal{H}_{\mathcal{K}})\simeq
H^2(C_{\underline{t}},\CC)\otimes_{\CC} \mathcal{K}$ it follows
that $\Phi_{\underline{t}}'(a,0,0,0,0)=a\Omega$. Hence it is
enough to prove that the induced map of quotient spaces
$\bar{\Phi'}_{\underline{t}}:\CC^4\to\mathbb{P}
H^2(C_{\underline{t}})$ is injective.

 The map $\bar{\Phi'}$ has a simple geometric description.
Indeed, fix a point $\underline{t^0}$ such that
$C_{\underline{t^0}}$ is smooth. Then for any $\underline{t}\in
\CC^4$ from a neighborhood of $\underline{t_0}$ we have a
$C^\infty$-diffeomorphism $\psi_{\underline{t}}$: $C_{\underline{t}}\to
C_{\underline{t_0}}$, which continuously depends on $\underline{t}$
and is the identity when $\underline{t}=\underline{t_0}$.
Let $\Omega_{\underline{t}}\in
H^2(C_{\underline{t}})$ be the 2-form given by formula
(\ref{2form}) (i.e. it is the form induced by the noncommutative
deformation of $\CC[C_{\underline{t}}]$). Hence we have a well
defined map $\Theta$: $\CC^4_{\underline{t}}\to
{\Bbb P}H^2(C_{\underline{t_0}})$:
$\Theta(\underline{t})=\psi_{\underline{t}}(\Omega_{\underline{t}})$
modulo scalars, and the map
$\bar{\Phi'}_{\underline{t}}$ is the differential of this map. It
is easy to see that we can extend the function $\Theta$ to
$\CC^4_{\underline{t}}\setminus\Sigma$ where $\Sigma$ is the set
of points $\underline{t}$ such that $C_{\underline{t}}$ is
singular, and this function is holomorphic in this domain (but
multivalued).

 To prove the theorem it is enough to show that the generic point
of $\CC^4_{\underline{t}}$ is not critical for $\Theta$. That is, it is
enough to show that $\Theta(\CC^4_{\underline{t}})$ cannot have
codimension $\ge 1$. Obviously if we show that the image of a
neighborhood of $\underline{t}=(1,1,1,1)$ is dense in a
neighborhood of $\Theta(\underline{t})$ we are done.

First let us notice that parameter $\underline{t}\in\CC^4$ gives a
bad parametrization on the space of the affine cubics with the
triangle at infinity. Indeed, the equation of the cubic
$C_{\underline{t}}$ has the form $$
\tilde{R}_p(X_1,X_2,X_3)=X_1X_2X_3+\sum_{i=1}^3(-X_i^2+p_i
(\underline{t})X_i)+p_0(\underline{t})+4,$$ where
$p_1=\bar{u}_0\bar{k}_0+\bar{k}_1\bar{u}_1$,
$p_2=\bar{u}_1\bar{u}_0+\bar{k}_0\bar{k}_1$,
$p_3=\bar{k}_0\bar{u}_1+ \bar{k}_1\bar{u}_0$, $p_0=
\bar{k}_0^2+\bar{k}_1^2+\bar{u}_0^2+\bar{u}_1^2-
\bar{k}_0\bar{k}_1\bar{u}_0\bar{u}_1$. Let us denote this cubic by
$\tilde{C}_p$. Unfortunately, the differential of the map
$\underline{t}\mapsto p(\underline{t})$ is zero at
$\underline{t}=(1,1,1,1)$, but obviously there exists a map
$\tilde{\Theta}$: $\CC_p^4\to \mathbb{P}H^2(C_{\underline{t_0}})$
such that $\Theta=\tilde{\Theta}\circ p$. It is enough to prove
that the image of a neighborhood of $p=0$ under the map
$\tilde{\Theta}$ is ``dense'', i.e. contains a ball.

The cubic surface $\tilde{C}_{0}$ is singular with four singular
points $S_\epsilon=2(\epsilon_1,\epsilon_2,\epsilon_3)$,
$\epsilon_i=\pm 1$, $\epsilon_1\epsilon_2\epsilon_3=1.$ It has a
simple geometric construction. The group $\ZZ_2$ acts on
$(\CC^*)^2$ by the inversion of both coordinates and there is an
obvious isomorphism between $(\CC^*)^2/\ZZ_2$ and $\tilde{C}_0$:
$z\mapsto (z_1^{-1}+z_1,z_2^{-1}+z_2, z_1z_2+z_1^{-1}z_2^{-1})$.
The four singular points $S_\epsilon$ are the images of four
$\ZZ_2$-fixed points $(\pm 1,\pm 1)\in(\CC^*)^2$. As
$(\CC^*)^2/\ZZ_2$ is homotopic to the two dimensional sphere, the
last description of $\tilde C_0$ implies that
$H_2(\tilde{C}_0,\CC)=\CC$.

Let $p\in \CC^4$ be small, and $\tilde{C}_p$ be smooth. The
homology group $H_2(\tilde{C}_p,\CC)$ is generated by the classes
of five spheres: by the classes of the four vanishing spheres
$\mathbb{S}^2_\epsilon$ and cycle $\mathbb{S}^2_0$. Here
$\mathbb{S}^2_\epsilon$ degenerates to the singular point
$S_\epsilon$ and $\mathbb{S}^2_0$ deforms to the generator of
$H_2(\tilde{C}_0)$ as $p$ tends to $0$.

As $\int_{\mathbb{S}^2_0}\Omega_p\to b\ne 0$ as $p\to 0$
(where $\Omega_p:=\Omega_{\underline{t}}$,
$p=p(\underline{t})$), we can write $$\tilde{\Theta}(p)=
(\int_{\mathbb{S}^2_{1,1,1}}
\Omega_p,\int_{\mathbb{S}^2_{1,-1,-1}}\Omega_p,
\int_{\mathbb{S}^2_{-1,1,-1}}\Omega_p,
\int_{\mathbb{S}^2_{-1,-1,1}}\Omega_p).$$
This formula shows that $\tilde{\Theta}$ has a
holomorphic extension at $p=0$.
We will now show that this point is not critical.

From corollary~\ref{sings} it follows that $p(\Sigma''_{1,1,1})$
is a union of four curves
$p(\Sigma''_{1,1,1})=\cup_{\epsilon\in\{\pm
1\}^3,\epsilon_1\epsilon_2\epsilon_3=1}K_{\epsilon}$,
\begin{equation}\label{param=s}
K_\epsilon=\{p,p_i=\epsilon_i s, i=1,2,3, p_0=s(2-s/2), s\in\CC\}.
\end{equation} The tangent vectors to $K_\epsilon$ at
the point $p=0$ form a basis in $\CC_p^4$. Hence it is enough to
prove that the derivatives of $\tilde{\Theta}$ along the curves
$K_\epsilon$ at $p=0$ are linearly independent.

The points of the curve $K_\epsilon$ correspond to the cubic
surfaces with three singular points which are close to
$S_{\epsilon'}$, $\epsilon'\ne \epsilon$. That is the homology
group $H_2(\tilde C_p,\CC)$ is two dimesional with generators
$[\mathbb{S}^2_0(p)]$ and $[\mathbb{S}^2_\epsilon(p)]$ and the
integral $\int_{\mathbb{S}_{\epsilon'}^2}\Omega_p$,
$\epsilon\ne\epsilon'$ equals identically zero along $K_\epsilon$.
Now let us study behavior of the integral
$I_\epsilon(p)=\int_{\mathbb{S}_\epsilon^2}\Omega_p$ as $p$ tends
to $0$ along the curve $K_\epsilon$.

We give an analysis only for $K=K_{1,1,1}$ because the other cases
are absolutely analogous. We will use the parametrization of
$K=\{p(s)\}_{s\in\CC}$ from formula (\ref{param=s}).

To estimate the integral $I_{(1,1,1)}(p)=I(p)$ we need to
understand the geometry of the universal family
$\mathcal{C}=\{(X,s)\in\CC^3\times\CC_s, R_{p(s)}(X)=0\}$ in the
neighborhood of $(X,s)=(2,2,2,0)$. After the change of the
variables
\begin{gather*}
X_i=\tilde{X_i}+2-s/2, \quad i=1,2,3,
\end{gather*}
the equation of the family $\mathcal{C}\subset\CC^3\times\CC_s$
takes the form:
\begin{multline*}
\tilde{X}_1\tilde{X}_2\tilde{X}_3+(2-s/2)(\tilde{X}_1\tilde{X}_2+
\tilde{X}_2\tilde{X}_3+\tilde{X}_3\tilde{X}_1)-\sum_{i=1}^3\tilde{X}^2_i
+ \\s^2/4\sum_{i=1}^3 \tilde{X}_i+s(2+s/4-s^3/8)=0.
\end{multline*}
Hence there are local coordinates $Y_i,s'$, $i=1,2,3$ at the
neighborhood of $(X,s)=(2,2,2,0)$ such that \begin{gather*}
X_i=Y_i+O(s)+O(|Y|^2), \quad |Y|^2=\sum_{i=1}|Y_i|^2,\\
s'=s+O(s^2)
\end{gather*} and the
equation of $\mathcal{C}$ at these coordinates takes the form
\begin{equation*}
Y_1^2+Y_2^2+Y_3^2-2(Y_1Y_2+Y_2Y_3+Y_3Y_1)=2s'.
\end{equation*}
For convenience we will make no difference between $s$ and $s'$.
After the linear change of the variables
\begin{multline*}
Z_1=\rm{i}(c^{-1}Y_2+c Y_3),\quad Z_2= \rm{i}(c^{-1}Y_3+c
Y_1),\quad Z_3=\rm{i}(c^{-1}Y_1+c Y_2),
\end{multline*}
with $c=(-1+\rm{i}\sqrt{3})/2$, we get the simplest possible
equation for $\mathcal{C}\cap \mathcal{U}$ (for a suitable
neighborhood $\mathcal{U}$ of zero):
\begin{equation*}
\sum_{i=1}^3 Z_i^2=2s.
\end{equation*}

 Having the description of $\mathcal{C}\cap\mathcal{U}$ from the
previous paragraph we can say that  the intersection of $\tilde
C_{p(s)}$, $s\ne 0$ with small neighborhood of $Z=0$ is
retractable to the two dimensional sphere $\mathbb{S}^2(s)= \{
Z\in \tilde C_{p(s)}|Im Z=0\}$ (see e.g. \cite{Arn}).
Obviously if $s>0$ then $\mathbb{S}^2(s)$ is the honest two
dimensional sphere of the radius $\sqrt{2s}$.

Let us write the form of $\Omega_{p(s)}$ in the coordinates $Z,s$.
As $$\frac{\partial\tilde{R}}{\partial X_3}=X_1X_2-2X_3+p_3(s)=
-2iZ_3+O(s^2)+O(|Z|^2)$$ we have $\Omega_{p(s)}=i dZ_1\wedge
dZ_2/2Z_3(1+O(s^2)+O(|Z|^2))$.

Suppose $s\in \mathbb{R}_+$. Now let us notice that if $Z\in
\mathbb{S}^2(s)$ then $|Z|^2\sim 2s$ and also form $dZ_1\wedge
dZ_2/Z_3$ is proportional to the standard volume form for the two
dimensional sphere. Hence we have
\begin{multline*}
I(s)=\int_{\mathbb{S}^2(s)}\Omega_{k(s)}=\int_{\mathbb{S}^2(s)} i
dZ_1\wedge dZ_2/2Z_3(1+O(s^2)+O(|Z|^2))=\\k \mathop{\rm
vol}(\mathbb{S}^2(s))(1+O(s)),
\end{multline*}
where $k\ne 0$. That is we have $I(s)=ks+O(s^2)$.

Thus we proved that $I'(s)\ne 0$.
The same is true for the other curves $K_\epsilon$.
Thus the map $\tilde{\Theta}$ is holomorphic at $p=0$ and this
point is not critical.
\end{proof}

\section{DAHA as universal deformation of $D_q\rtimes\CC[\ZZ_2]$}
The algebra $D_q$ has the generators $X^{\pm 1},P^{\pm 1}$ and
these elements satisfy the defining relation: $$PX=qXP.$$ Let us
fix notation $s$ for the generator of $\ZZ_2$ and $e$ for the
unit. We use the notation $\tilde{D_q}$ for
$D_q\rtimes\CC[\ZZ_2]$.

First we will calculate the homology $HH_*(\tilde{D_q})=
H_*(\tilde{D_q},\tilde{D_q})$. We calculate homology instead of
cohomology just to shorten notations. It is easy to see that  the
same method works for cohomology and that
$HH_j(\tilde{D_q})=HH^{2-j}(\tilde{D_q})$.

\begin{prop}\label{HHDq} If $q$ is not a root of unity we have:
\begin{gather*}
HH_0(\tilde{D_q})=HH^2(\tilde{D_q})=\CC^5,\\
HH_1(\tilde{D_q})= HH^1(\tilde{D_q})=0,\\
HH_2(\tilde{D_q})=HH^0(\tilde{D_q})=\CC.
\end{gather*}
\end{prop}

\subsection{Spectral sequence}\label{spseq}
There is a natural structure of a $\ZZ_2$-module on the homology
$HH_i(D_q,gD_q )$, where $g=e,s$ is one of the elements of
$\ZZ_2$. More precisely there is an action of $\ZZ_2$ on the
standard Hochschild complex for $HH_i(D_q,gD_q)$ by the formulas:
$$g\cdot(m\otimes a_1\otimes\dots \otimes a_r)=m^{g}\otimes
a_1^g\otimes\dots\otimes a_r^g.$$

 Proposition 3.1 from the paper \cite{Al} implies:
\begin{prop}There is a decomposition:
$$HH_*(\tilde{D_q},\tilde{D_q})=
HH_*(\tilde{D_q},\tilde{D_q})_{[e]}\oplus HH_*(\tilde{D_q},
\tilde{D_q})_{[s]},$$ and the spectral sequence: $$
E^2_{r,s,[g]}=H_r(\ZZ_2, H_s(D_q, gD_q))\Rightarrow
HH_{r+s}(\tilde{D_q},\tilde{D_q})_{[g]}$$ The same statement holds
for the cohomology.
\end{prop}
For calculation of $H_*(D_q,gD_q)$ we will use the Koszul
resolution.
\begin{rema}The  elementary group theory implies
$$H_{>0}(\ZZ_2,V)=0, H_0(\ZZ_2,V)=V^{\ZZ_2}$$ for any
$\ZZ_2$-module $V$ (over $\CC$).
\end{rema}

\subsection{Resolution}
 Let us denote by $D^e_q$ algebra $D_q\otimes D_q^{opp}$, where
$D_q^{opp}$ is the algebra $D_q$ with the opposite multiplication.
The elements $p=P\otimes P^{-1}-1$, $x=X\otimes X^{-1}-1$ commute
and $D_q^e/{\bf I}=D_q$, where ${\bf I}=(x,p)$ is the
$D_q^e$-submodule generated by these elements. Hence the
corresponding Koszul complex yields a free resolution $W_*$ of
$D^e_q$-module $D_q$:
$$D^e_q\stackrel{d_1}{\to} D^e_q\oplus
D^e_q\stackrel{d_0}{\to}D^e_q\stackrel{\mu}{\to} D_q,$$ where
$\mu(X^iP^j\otimes P^{j'}X^{i'})=X^i P^{j+j'}X^{i'}$ and for
$z=z_1\otimes z_2$ we have $d_0(z,0)=zp=z_1 P\otimes P^{-1}
z_2-z_1\otimes z_2$, $d_0(0,z)=zx=z_1X\otimes X^{-1}z_2-z_1\otimes
z_2$, $d_1(z)= (zx,-zp)=(z_1X\otimes X^{-1}z_2-z_1\otimes
z_2,-z_1P\otimes P^{-1} z_2+z_1\otimes z_2).$

\subsection{Calculation of $HH_*(D_q,D_q)$}
After  multiplication of the resolution $W_*$ by the
$D^e_q$-module $D_q$ we get the complex of $D_q$ modules:
$$D_q\otimes_{D_q^e} W_*: 0\to D_q\stackrel{\hat{d}_1}{\to}
D_q\oplus D_q\stackrel{\hat{d}_0}{\to}D_q\to 0,$$
$\hat{d}_1(z)=(XzX^{-1}-z,-PzP^{-1}+z)$, $\hat{d}_0(z_1,z_2)=
Pz_1P^{-1}-z_1+Xz_2X^{-1}-z_2$. The homology of this complex
yields the Hochschild homology: $HH_i(D_q)=H_i(D_q\otimes_{D_q^e}
W_*)$.

We have $H_2(D_q\otimes_{D_q^e} W_*)=\ker \hat{d}_1=\CC$.

The element $(z_1,z_2)$, $z_k=\sum c_k^{ij}X^iP^j$ is from
$\ker\hat{d}_0$ if and only if $$
c_1^{ij}(1-q^i)+c_2^{ij}(1-q^j)=0$$ for all $(i,j)\in \ZZ^2$. As
the image $Im \hat{d}_1$ is spanned by the vectors $$\sum c^{ij}
X^iP^j (q^j-1,1-q^i)$$ we get $H_1(D_q\otimes_{D_q^e}
W_*)=\CC^2=\langle (1,0)+\hat{d}_1(D_q),
(0,1)+\hat{d}_1(D_q)\rangle$.

It is easy to see that $Im\hat{d}_0=\{\sum c^{ij}
X^iP^j|c^{00}=0\}$. Hence $$H_0(D\otimes_{D_q^e} W_*)=\CC=\langle
1+\hat{d}_0(D_q\oplus D_q)\rangle.$$

\subsection{Calculation of $HH_*(D_q,sD_q)$}
After  multiplication of the resolution $W_*$ by the
$D^e_q$-module $sD_q$ we get the complex of $D_q$ modules:
$$sD_q\otimes_{D_q^e} W_*: 0\to D_q\stackrel{\bar{d}_1}{\to}
D_q\oplus D_q\stackrel{\bar{d}_0}{\to}D_q\to 0,$$
$\bar{d}_1(z)=(X^{-1}zX^{-1}-z,-P^{-1}zP^{-1}+z)$,
$\bar{d}_0(z_1,z_2)= P^{-1}z_1P^{-1}-z_1+X^{-1}z_2X^{-1}-z_2$. The
homology of this complex yield the Hochschild homology:
$H_i(D_q,sD_q)=H_i(sD_q \otimes_{D_q^e} W_*)$.

Let  $F=\{ c\in Fun_{fin}(\ZZ^2,\CC)\}$ be the space of the
functions with the finite support. Let us introduce two
"differentiations" on this space
$$(\delta^{(1)}c)(k,l)=q^{-l}c(k+2,l)-c(k,l),\quad
(\delta^{(2)}c)(k,l)=q^{k}c(k,l+2)-c(k,l).$$ A simple calculation
shows that
\begin{gather*} Im\delta^{(1)}=\{c\in
F|\sum_{i=-\infty}^\infty
q^{-2ij}c(\varepsilon_1+2i,\varepsilon_2+2j)=0,
\forall j\in \ZZ, \varepsilon_1,\varepsilon_2=0,1\},\\
Im\delta^{(2)}=\{c\in F|\sum_{j=-\infty}^\infty q^{-2ij}
c(\varepsilon_1+2i,\varepsilon_2+2j)=0,\forall i\in
\ZZ,\varepsilon_1,\varepsilon_2=0,1\}.
\end{gather*}

 Using these
operation we can rewrite the formulas for $\bar{d}_i$:
\begin{gather*}
\bar{d}_1(\sum c(k,l)X^k P^l)=(\sum(\delta^{(1)}c)(k,l)
X^kP^l,-\sum(\delta^{(2)} c)(k,l)X^kP^l),\\ \bar{d}_0(\sum
c_1(k,l)X^kP^l,\sum c_2(k,l) X^kP^l)=\sum
q^{-2ij}(\delta^{(2)}c_1+\delta^{(1)}c_2)(k,l)X^kP^l.
\end{gather*}

Obviously we have $\ker\bar{d}_1=H_2(D_q,sD_q)=0$.

The kernel of $\bar{d}_0$ consists of pairs $(c_1,c_2)\in F\oplus
F$ satisfying the equation $\delta^{(2)}c_1=-\delta^{(1)}c_2$.
Hence $(c_1,c_2)\in \ker \bar{d}_0$ implies
$$\delta^{(2)}
(\sum_{i=-\infty}^\infty q^{-2ij}
c_1(\varepsilon_1+2i,\varepsilon_2+2j))= -\sum_{i=-\infty}^\infty
q^{-2ij}(\delta^{(1)} c_2)(\varepsilon_1+2i,\varepsilon_2+2j)=0,$$
for all $j\in\ZZ, \varepsilon_1,\varepsilon_2=0,1$. As $c_1$ has
the finite support it implies $$\sum_{i=-\infty}^\infty
q^{-2ij}c_1(\varepsilon_1+2i,\varepsilon_2+2j)=0,$$ i.e.
$c_1=\delta^{(1)}c$ for some $c\in F$. Hence
$\delta^{(1)}\delta^{(2)}c=\delta^{(2)}\delta^{(1)}c=-\delta^{(1)}c_2$
because $\delta^{(1)}$ and $\delta^{(2)}$ obviously commute. For
the same reason as before $\delta^{(2)}c=-c_2$. Thus we proved
that $Ker \bar{d}_0= Im\bar{d}_1$ and $H_1(D_q,sD_q)=0$.

We have $Im\bar{d}_0=Im\delta^{(1)}+Im\delta^{(2)}$ where the sum
is not direct. That is $$ Im\bar{d}_0=\{\sum  c(k,l)
X^kP^l|\sum_{i,j\in\ZZ}  q^{-2ij}c(\varepsilon_1+2i,
\varepsilon_2+2j)=0,\forall \varepsilon_1,\varepsilon_2=0,1\}$$ and
$H_0(D_q,sD_q)=\CC^4$ is spanned by four classes
$X^{\varepsilon_1}P^{\varepsilon_2}+\bar{d}_0(D_q\oplus D_q)$,
$\varepsilon_i=0,1$.

\subsection{} To complete the calculation of the $E_2$ term of the
spectral sequence we need to describe the action of $\ZZ_2$ on the
homology $H_i(D_q,gD_q)$.

Let us mention that the Koszul resolution:
$$W_2\stackrel{d_1}{\to} W_1\stackrel{d_0}{\to} W_0\to D_q,$$
is actually a resolution in the abelian category of the
$D_q^e$-modules with $\ZZ_2$-action, where $\ZZ_2$-action on $W_i$
is given by the formulas
\begin{gather*}
z\in W_2,\quad s\cdot z=z^s(X^{-1}P^{-1}\otimes PX),\\
(z_1,z_2)\in W_1,\quad s\cdot (z_1,z_2)=-(z_1^s (P^{-1}\otimes P),
z_2^s(X^{-1}\otimes X)),\\
z\in W_0,\quad s\cdot z=z^s.
\end{gather*}
By the way this resolution extends the standard $\ZZ_2$-action on
$D_q$.

The action of $\ZZ_2$ on the standard Hochschild complex from the
 subsection~\ref{spseq} is the action of $\ZZ_2$ on the
 derived functor of
tensor multiplication in the category of $D_q^e$-modules with
$\ZZ_2$-action and on $HH_0(D_q,gD_q)$ this action is induced by
the standard one. As the derived functor does not depend on the
used resolution the corresponding action of $\ZZ_2$ on the complex
$W'_*=W_*\otimes_{A^e} gD_q$ induces the desired action on
$HH_i(D_q,gD_q)$. This $\ZZ_2$-action is given by the formulas:
\begin{gather*}
z\in W'_2,\quad s\cdot z=P^gX^gz^sX^{-1}P^{-1},\\ (z_1,z_2)\in
W'_1,\quad s\cdot (z_1,z_2)=-(P^gz^sP^{-1},X^gz^sX^{-1}),\\ z\in
W'_0,\quad s\cdot z=z^{s}.
\end{gather*}

 From the results of the previous two subsections we see that the
action of $\ZZ_2$ on $H_0(D_q,D_q)$, $H_2(D_q,D_q)$,
$H_0(D_q,sD_q)$ is trivial and $\ZZ_2$ acts on $H_1(D_q,D_q)$ by
multiplication by $-1$. Hence we get
\begin{gather*}
H_{>0}(\ZZ_2,H_*(D_q,gD_q))=0,\\
H_0(\ZZ_2,H_0(D_q,D_q))=H_0(\ZZ_2,H_2(D_q,D_q))=\CC,\\
H_0(\ZZ_2,H_1(D_q,D_q))=H_0(\ZZ_2,H_1(D_q,sD_q))=0,\\
H_0(\ZZ_2,H_2(D_q,sD_q))=\CC^4.
\end{gather*}

This calculation completes the proof of Proposition~\ref{HHDq}.

\subsection{Universal property of DAHA}
In this subsection we prove
\begin{thm}\label{univq} If  $q$ is not a root of unity then the
 family $\{
H(\underline{t};q)\}_{\underline{t}\in\CC^4,q\in\CC}$ gives a
universal deformation of the algebra $\tilde{D}_q=\CC_q[X^{\pm
1},P^{\pm 1}]\rtimes \CC[\ZZ_2]$, and the family $\{
eH(\underline{t};q)e\}_{\underline{t}\in\CC^4,q\in\CC}$ gives a
universal deformation of the algebra $\CC_q[X^{\pm 1},P^{\pm
1}]^{\ZZ_2}$.
\end{thm}

 Let us remind the definition of a {\it universal deformation} .
The flat $R$-algebra $A_R$ (with $R$ being a local commutative
Artinian algebra and $\mathfrak{m}\subset R$ is the maximal ideal)
together with an isomorphism $A_R/\mathfrak{m}\simeq A$ is called
a deformation of $A$ over $S=Spec(R)$.  $A_R$  is a universal
deformation of $A$ if for every deformation $A_{\mathcal O(S)}$ of
$A$ over an Artinian base $S$ there exists a map $\tau$: $S\to
Spec(R)$ such that isomorphism $A\simeq A_R/\mathfrak{m}$ extend
to isomorphism $A_S\simeq \tau^* A_R$.

Let  $v=(\underline{t'},q')$ be a nonzero vector and
$H'=H(\underline{1}+\epsilon\underline{t}';q+\epsilon
q')/\epsilon\cdot
H(\underline{1}+\epsilon\underline{t}';q+\epsilon q')$ is the
$\CC[\epsilon]/(\epsilon^2)$-algebra. The theorem basically
follows from the calculation of the Hochschild cohomology and the
following lemma

\begin{lem}\label{nontriv}
 There is no isomorphism of $\CC[\epsilon]/(\epsilon^2)$-algebras
between $H'$ and
$H(\underline{1};q)\otimes_{\CC}\CC[\epsilon]/(\epsilon^2)$ which
is  equal to the identity map modulo the ideal $(\epsilon)$.
\end{lem}

In the proof of this lemma it is more convenient to use the
following description of the algebra $H(\underline{t};q)$.

\begin{prop}
The algebra $H(\underline{t};q)$ is generated by elements
$Y=V_1V_0$, $T=V_1$, $X=q^{1/2}V_0V_0^{\vee}$, modulo the defining
relations
\begin{gather*}
XT=T^{-1}X^{-1}+(u_1^{-1}-u_1),\\
Y^{-1}T=T^{-1}Y+(k_0^{-1}-k_0),\\
(T-k_1)(T+k_1^{-1})=0,\\
YX=qT^2XY+q(u_1-u_1^{-1})TY+(k_0-k_0^{-1})TX+q^{1/2}(u_0-u_0^{-1})T.
\end{gather*}
\end{prop}

 This description is more convenient because the generators
 $X,Y,T$ tend to  $X,P,s$ as $\underline{t}$ tends to $\underline{1}$.

\begin{proof}[Proof of the Lemma~\ref{nontriv}]
Let us denote by $\phi$ the natural isomorphism of vector spaces
$$H(\underline{1};q)\to
\epsilon\cdot H(\underline{1}+\epsilon\underline{t}';q+\epsilon q'
)/\epsilon^2 \cdot
H(\underline{1}+\epsilon\underline{t};q+\epsilon q').$$ Assume
that there is an isomorphism between $H'$ and
$H(\underline{1};q)\otimes_{\CC}\CC[\epsilon]/(\epsilon^2)$
lifting the identity. Hence there exists a linear map $f$:
$\CC_q[X^{\pm 1},P^{\pm 1}]\to \CC_q[X^{\pm 1},P^{\pm
1}]\rtimes\CC[\ZZ_2]$ such that the equation
\begin{multline*}
(P+\phi\circ f(P))(X+\phi\circ f(P))=(q+\epsilon q')(1+2\epsilon
k_1's)(X+\phi\circ f(X))\times\\
(P+\phi\circ f(P))+ 2\epsilon q u'_1sP+2\epsilon k'_0 sX+2\epsilon
q^{1/2} u'_0s,
\end{multline*}
holds modulo $\epsilon^2$. Taking the first order term in
$\epsilon$ and applying $\phi^{-1}$, we get
\begin{multline*}
(f(P)X-qXf(P))+(Pf(X)-qf(X)P)=q'XP+2k'qX^{-1}P^{-1}s+\\2q
u'_1P^{-1}s +2k'_0X^{-1}s +2q^{1/2}u'_0 s.
\end{multline*}

Now let us show that the last equation has no solutions unless
$(\underline{t}',q')=0$. For that let us introduce the functionals
$I,I_{\delta_1,\delta_2}$, $\delta_i=0,1$ on $\CC_q[X^{\pm
1},P^{\pm 1}]\rtimes \CC[\ZZ_2]$:
\begin{gather*}
I(h)=c^0_{1,1},\quad I_{\delta_1,\delta_2}(h)=\sum_{i,j\in\ZZ}
q^{i+j-2ij}c^1_{\delta_1+2i,\delta_2+2j},
\end{gather*}
where $h=\sum_{i,j\in \ZZ,\epsilon=0,1}c_{ij}^\epsilon
X^iP^js^{\epsilon}$. It is easy to see that application of
$I,I_{\delta_1,\delta_2}$ to the LHS of the last equation yields
zero. Indeed let us check that
$I(\tilde{f})=I_{\delta_1,\delta_2}(\tilde{f})=0$,
($\delta_1,\delta_2=0,1$) for $\tilde{f}= Pf(X)-qf(X)P$ (the
calculation for $f(P)X-qXf(P)$ is absolutely analogous). Obviously
it is enough to check it in the case $f=X^iP^js^\epsilon$, $i,j\in
\ZZ$, $\epsilon=0,1$.

First let us consider the case when $\epsilon=0$. In this case
$I_{\delta_1,\delta_2}=0$ for obvious reasons and $I(\tilde{f})=0$
because $ \tilde{f}=(q^i-q)X^iP^{j+1}.$ In the case $\epsilon=1$
we have $I(\tilde{f})=0$  for obvious reasons and
$I_{\delta_1,\delta_2}(\tilde{f})=0$ because $\tilde{f}=
q^iX^iP^{j+1}s-qX^iP^{j-1}s$.

 At the same time application of the functionals to the RHS
is zero if and only if $(\underline{t}',q')=0$.
\end{proof}
\begin{proof}[Proof of the theorem~\ref{univq}] We know that
$HH^3(\tilde{D}_q,\tilde{D}_q)=0$, hence all deformations of
$\tilde{D}_q$ are unobstructed. Thus the second Hochschild
cohomology $HH^2(\tilde{D}_q,\tilde{D}_q)=\CC^5$ is the tangent
space to the moduli space of all deformations. The deformations
coming from the family
$\{H(\underline{t},q)\}_{\underline{t}\in\CC^4,q\in\CC}$ yield a
subspace in the moduli space of all deformations. The last lemma
shows that this subspace is of dimension $5$.

The second part of the statement follows from the existence of the
isomorphism $e\tilde{D}_qe\simeq \CC_q[X^{\pm 1},P^{\pm
1}]^{\ZZ_2}$. Indeed this isomorphism implies that the family
$\{eH(\underline{t};q)e\}_{\underline{t}\in\CC^4,q\in \CC}$ gives
a flat deformation family of the algebra $\CC_q[X^{\pm 1},P^{\pm
1}]^{\ZZ_2}$. The rest of the proof is absolutely the same because
the Morita equivalence of $\tilde{D}_q$ and $\CC_q[X^{\pm
1},P^{\pm 1}]^{\ZZ_2}$ implies that
$HH^*(\tilde{D}_q)=HH^*(\CC_q[X^{\pm 1},P^{\pm 1}]^{\ZZ_2})$.
\end{proof}

\section{Poisson deformations of $\CC[X^{\pm 1},P^{\pm 1}]^{\ZZ_2}$}
In this section we prove that the family of structure rings of the
cubic surfaces $\tilde{C}_p$ equipped with the Poisson two-form
$r\Theta=r\Omega^{-1}$, where $r\in\CC$ and $\Omega$ is given by
(\ref{2form}), yields a universal formal Poisson deformation of
the ring $\CC[X^{\pm 1},P^{\pm 1}]^{\ZZ_2}$ with the  Poisson
structure induced by the Poisson structure $$ \{ X,P\}=XP,$$ on
$\CC[X^{\pm 1},P^{\pm 1}]$. Our arguments are analogous to the
arguments from the third section of the paper \cite{GK}.
Particularly, the proof of the lemma \ref{HP} is the parallel to
the proof of the lemma 3.1 from the paper \cite{GK}.

Let us remind some standard definitions from  deformation theory.
In this sections all rings and algebras are commutative. Let $R$
be an algebra. Recall that $A$ is said to be a {\it Poisson
$R$-algebra}, if $A$ is equipped with an $R$-linear skew-symmetric
bracket $\{\cdot,\cdot\}$ that satisfies the Leibniz rule and the
Jacobi identity. Let us denote by $\mathcal{O}_{p,r}$,
$p\in\CC^4_p$, $r\in \CC$ the structure ring of the cubic surface
$\tilde{C}_p$ equipped with the Poisson structure induced by the
form $r\Omega$, where $\Omega$ is given by the formula
(\ref{2form}).

\begin{defi} A {\it Poisson deformation} of a Poisson algebra $A$
over the spectrum $S=Spec R$ of a local Artinian algebra $R$ with
maximal ideal $\mathfrak{m}\subset R$ is a pair of a flat Poisson
$R$-algebra $A_R$ and a Poisson $R$-algebra isomorphism
$A_R/\mathfrak{m}\simeq A$.
\end{defi}

A Poisson $R$-algebra $A_R$ over a complete local $\CC$-algebra
$\langle R,\mathfrak{m}\rangle$  is called a {\it universal formal
Poisson deformation} of the Poisson algebra $A$ if for every
Poisson deformation $A_{\mathcal{O}(S)}$ over a local Artinian
base $S$ there exists a unique map $\tau$: $S\to Spec R$ such that
the isomorphism $A\simeq A_R/\mathfrak{m}$ extends to a Poisson
isomorphism:
$$ A_{\mathcal{O(S)}}\simeq\tau^*A_R.$$

\begin{rema} If we rewrite the previous definitions without
mentioning the Poisson structure we get the definitions of a
formal deformation and universal formal deformation for an affine
scheme.
\end{rema}

We can consider the family $\mathcal{O}_{p,r}$,
$p\in\CC^4,r\in\CC$ of algebras as  a $\CC[p,r]$-Poisson algebra
$\mathcal{O}$. Let us denote by $\hat{\mathcal{O}}$ the completion
of this ring with respect to the variables $p,r$.

\begin{thm} The Poisson $\CC[[p,r]]$-algebra $\hat{\mathcal{O}}$ is a
universal Poisson deformation of the Poisson algebra $\CC[X^{\pm
1},P^{\pm 1}]^{\ZZ_2}$.
\end{thm}

We will prove this theorem using the technique of  Poisson
cohomology $HP^{\bullet}(A)$ (here $A$ is a Poisson algebra). The
reader can find the definition and main properties of  Poisson
cohomology in the Appendix to paper \cite{GK}.

Basically the Poisson cohomology $HP^{\bullet}(A)$ is the total
cohomology of the bicomplex
$$ DP^{l,k}(A)=Hom_A(\Lambda_A^k Har_l(A),A).$$
In this bicomplex  $Har_{\bullet}(A)$, $d: Har_{\bullet}(A)\to
Har_{\bullet+1}$ is the Harrison complex which is quasiisomorphic
to the co-called  cotangent complex $\Omega_{\bullet}(A)$ in the
case when $Spec(A)$ is  smooth (i.e. the complex
$Har_{\bullet}(A)$ has nonzero cohomology only at zero degree and
this cohomology are equal to $\Omega^1_{Spec(A)}$). For the
definition of the complex $Har_{\bullet}(A)$ see the Appendix of
the paper \cite{GK}; for our purposes it is enough to know only
the first two terms of the complex $Har_{\bullet}(A)$:
$$\longrightarrow S^2A\otimes A\stackrel{d}{\longrightarrow}
A\otimes A\longrightarrow 0, $$ with the differential given by
$d$: $(a\otimes b+b\otimes a) \otimes c\mapsto a\otimes
bc+b\otimes ac-ab\otimes c.$ By the way, from this description we
see that the zero homology of the complex $Har_{\bullet}$ is the
module $\Omega^1A$ of K\"ahler differentials.

The Poisson structure $\Theta\in DP^{0,2}(A)\simeq
Hom_{\CC}(\Lambda^2 A,A)$ yields the second differential $\delta$:
$DP^{\bullet,\bullet}(A)\to DP^{\bullet,\bullet+1}(A),$
$\delta(a)=\{\Theta,a\},$ where $\{\cdot,\cdot\}$ is Gerstenhaber
bracket.

The most important property of  Poisson cohomology is the fact
that the second Poisson cohomology $HP^2(A)$ controls the first
order formal Poisson deformations of the algebra $A$:

\begin{thm}\cite{GK}
Let $A$ be a Poisson algebra. Assume that $HP^1(A)=0$ and that
$HP^2(A)$ is a finite-dimensional vector space over $\CC$. Then
there exists a closed subscheme $S\subset \widehat{HP^2(A)}$ and a
Poisson $\CC[S]$ algebra $A_S$ which is a universal formal Poisson
deformation of the algebra $A$.
\end{thm}

The theorem implies that there is a map $$ \varphi:
\widehat{\CC^5_{p,t}}\to S\subset \widehat{HP^2(A)},$$ and we only
need to show that $\varphi(\widehat{\CC^5_{p,t}})=S$ and $\varphi$
is injective. More exactly, we will prove

\begin{lem}\label{HP} We have $HP^1(\CC[X^{\pm 1},P^{\pm 1}]^{\ZZ_2})=0$,
$HP^2(\CC[X^{\pm 1},P^{\pm 1}]^{\ZZ_2})=\CC^5$ and the map
$\varphi$ is an isomorphism.
\end{lem}

Before proving the lemma let us remark that the family
$\{C_p\}_{p\in \CC^4}$ yields the deformation of the affine scheme
$Z=Spec(\CC[X^{\pm 1},P^{\pm 1}]^{\ZZ_2})$. The tangent space to
the moduli space of the formal deformations of the affine scheme
$Z$ is equal to $Ext^1(\Omega^1_Z,\mathcal{O}_Z)$. The standard
calculations from deformation theory (see for example exercises
after the chapter 16 in the book \cite{E}) imply

\begin{lem}\label{algdef} The natural map: $\widehat{\CC_p^4}\to
Ext^1(\Omega^1_Z,\mathcal{O}_Z)$ is an isomorphism.
\end{lem}

\begin{cor} The family $\{ C_p\}_{p\in \CC^4_p}$ form a universal
formal deformation of the scheme $Z$.
\end{cor}

\begin{proof}[Proof of the lemma~\ref{HP}] Let us denote by $A$
algebra $\CC[X^{\pm 1},P^{\pm 1}]^{\ZZ_2}$. First let us remark
that the calculation  inside the proof of  theorem \ref{maindef}
(more precisely the fact that the point $0$ is not critical for
the map $\tilde{\Theta}$) implies that the map $\varphi$:
$\widehat{\CC^5_{p,r}}\to \widehat{HP^2(A)}$ is an embedding.
Hence $\dim HP^2(A)\ge 5.$ Now we prove that $\dim HP^2(A)\le 5$.

 Immediately from the definition of the bicomplex
$DP^{\bullet,\bullet}$ we get:
$$ DP^{0,0}(A)\simeq A,\quad DP^{0,>0}(A)=0,\quad
DP^{1,\bullet}(A)\simeq RHom^{\bullet}(\Omega^1A,A).$$


The first term of the spectral sequence $E^{\bullet,\bullet}$
associated to the bicomplex $HP^{\bullet,\bullet}(A)$ is of the
form:
$$ E_1^{p,q}=Ext^q(\Lambda^p\Omega^1A,A),$$
because $Spec(A)$ is a complete intersection and hence the complex
$Har_\bullet(A)$ has nontrivial cohomology only at the degree zero
which is equal to the module $\Omega^1A$ of K\"ahler
differentials.

Thus $E_1^{0,1}=Ext^1(A,A)=0$ and $E^{0,2}_1=Ext^2(A,A)=0$ and by
Lemma~\ref{algdef} $E_1^{1,1}=\CC^4$. Let us calculate
$E_1^{p,0}$.

Let $U$ be the subset of the smooth points of $Z=Spec(A)$ and $j$:
$U\to Z$ be the corresponding embedding. Then we have:
\begin{multline*}
E^{p,0}_1\simeq Hom(\Lambda^p\Omega^1_Z,\mathcal{O}_Z)\simeq
Hom(\Lambda^p\Omega^1_Z,j_*\mathcal{O}_Z)\simeq\\
Hom_U(\Lambda^p\Omega_U,\mathcal{O}_U)\simeq
H^0(U,\Lambda^p\mathcal{T}_U),
\end{multline*}
where $\mathcal{T}_U$ is tangent bundle to the scheme $U$. As
quotient the map $\pi$: $(\CC^*)^2\to Z$ is \'etale  over $U$ we
have
$$H^0(U,\Lambda^p\mathcal{T}_U)\simeq
H^0(\pi^{-1}(U),\Lambda^p\mathcal{T}_{(\CC^*)^2})^{\ZZ_2}.$$

Now we are ready to calculate $E_2$. The complement
$(\CC^*)^2\setminus U$ has codimension $2$ and differential $d$:
$E_1^{p,0}\to E_1^{p+1,0}$ is induced by the Poisson differential
on $H^0((\CC^*)^2,\Lambda^p\mathcal{T}((\CC^*)^2))$, hence:
$$ E_2^{p,0}\simeq HP^p((\CC^*)^2)^{\ZZ_2}\simeq
H^p((\CC^*)^2)^{\ZZ_2}.$$ That is we have $E_2^{1,0}=E_2^{0,1}=0$
and $E^{2,0}_2=\CC$, $E^{0,2}_2=0$, $\dim E_2^{1,1}\le 4$. Thus we
get $HP^1(A)=0$, $\dim HP^2(A)\le 5$.
\end{proof}

\section{Appendix}
\subsection{Proof of Theorem \ref{eqnCub}}
Let us prove the cubic relation between the generators
$X_1,X_2,X_3$ of the center of the DAHA.

As elements $X_i$ are central, we have
\begin{gather}\label{X1per}
X_1=(V_1^{\vee})^{-1}X_1V_1^{\vee}=V_1V_1^{\vee}+(V_1^{\vee})^{-1}V_1^{-1},\\
\label{X3per} X_3=V_1^{-1}X_3V_1=V_0^{\vee}V_1+
V_1^{-1}(V_0^{\vee})^{-1}.
\end{gather}

\begin{multline}\label{X1X2}
X_1X_2=(V_1V_1^{\vee}+(V_1^{\vee})^{-1}V_1^{-1})
(V_1^{\vee})^{-1}(V_0^{\vee})^{-1}+V_0^{\vee}V_1^{\vee}(V_1V_1^{\vee}+\\
(V^{\vee}_1)^{-1}V_1^{-1})=V_1(V_0^{\vee})^{-1}+V_0^{\vee}V_1^{-1}+
(V^{\vee}_1)^{-1}V_0+V_0^{-1}V_1^{\vee}\\=X_3-\bar{u}_0\bar{k}_1+
(V^{\vee}_1)^{-1}V_0+V_0^{-1}V_1^{\vee};
\end{multline}
here the first equality uses (\ref{X1per}) and the fact that
$X_1\in Z$; the second uses (\ref{VVVV}) twice; the third uses
(\ref{u0}), (\ref{k1}).

Our strategy is to extract the sum $X_1^2+X_2^2+X_3^2$ from
$X_1X_2X_3$ and to see what remains. Following this strategy we
expand:
\begin{multline*}
(V^{\vee}_1)^{-1}V_0X_3=(V_1^{\vee})^{-1}V_0(V_0^{\vee}V_1+
V_1^{-1}(V_0^{\vee})^{-1})=\\(V_1^{\vee})^{-1}V_0V_0^{\vee}V_1+
(V_1^{\vee})^{-1}V_0V_1^{-1}(V_0^{\vee})^{-1},
\end{multline*}
here we use (\ref{X3per}) in the first equality. Now let us expand
two summands separately:
\begin{multline*}
(V_1^{\vee})^{-1}V_0V_0^{\vee}V_1=(V_1^{\vee})^{-1}
V_1^{-1}(V_1^{\vee})^{-1}V_1=
(V_1V_1^{\vee})^{-2}+\\\bar{k}_1(V_1^{\vee})^{-1}V_1^{-1}(V_1^{\vee})^{-1},
\end{multline*}
here the first equality follows from (\ref{VVVV}) second from
(\ref{k1});
\begin{multline*}
(V_1^{\vee})^{-1}V_0V_1^{-1}(V_0^{\vee})^{-1}=
(V_1^{\vee})^{-1}V_0^{-1}V_1^{-1}(V_0^{\vee})^{-1}+
\bar{k}_0(V_1^{\vee})^{-1}V_1^{-1}(V_0^{\vee})^{-1}=\\
(V_1^{\vee})^{-1}V_0^{\vee}V_1^{\vee}(V_0^{\vee})^{-1}+
\bar{k}_0(V_1^{\vee})^{-1}V_1^{-1}(V_0^{\vee})^{-1}=
(V_0^{\vee}V_1^{\vee})^{-2}+\bar{u}_1(V_1^{\vee})^{-1}(V_0^{\vee})^{-2}+\\
\bar{u}_0(V_1^{\vee})^{-2}(V_0^{\vee})^{-1}+
\bar{u}_0\bar{u}_1(V_1^{\vee})^{-1}(V_0^{\vee})^{-1}+
\bar{k}_0(V_1^{\vee})^{-1}V_1^{-1}(V_0^{\vee})^{-1}=
(V_0^{\vee}V_1^{\vee})^{-2}-\\
\bar{u}_0\bar{u}_1(V_1^{\vee})^{-1}(V_0^{\vee})^{-1}+
\bar{u}_1(V_1^{\vee})^{-1}+\bar{u}_0(V_0^{\vee})^{-1}+
\bar{k}_0(V_1^{\vee})^{-1}V_1^{-1}(V_0^{\vee})^{-1}
\end{multline*}
the first equality is (\ref{k0}), the second is (\ref{VVVV}), the
third and fourth are (\ref{u0}), (\ref{u1}) (twice).

Now do the same procedure with $V_0^{-1}V_1^{\vee}X_3=
X_3V_0^{-1}V_1^{\vee}$:
\begin{multline*}
X_3V_0^{-1}V_1^{\vee}=(V_0^{\vee}V_1+V_1^{-1}(V_0^{\vee})^{-1})
V_0^{-1}V_1^{\vee}=
V_0^{\vee}V_1V_0^{-1}V_1^{\vee}+\\V_1^{-1}(V_0^{\vee})^{-1}
V_0^{-1}V_1^{\vee},
\end{multline*}
here the first equality is (\ref{X3per}). Expanding the two
summands separately we get:
\begin{multline*}
V_0^{\vee}V_1V_0^{-1}V_1^{\vee}=V_0^{\vee}V_1V_0V_1^{\vee}-
\bar{k}_0V_0^{\vee}V_1V_1^{\vee}=
V_0^{\vee}(V_1^{\vee})^{-1}(V_0^{\vee})^{-1}V_1^{\vee}-\\
\bar{k}_0V_0^{\vee}V_1V_1^{\vee}=(V_0^{\vee}V_1^{\vee})^2-
\bar{u}_1(V_0^{\vee})^2V_1^{\vee}-\bar{u}_0V_0^{\vee}(V_1^{\vee})^2+
\bar{u}_0\bar{u}_1V_0^{\vee}V_1^{\vee}-\\
\bar{k}_0V_0^{\vee}V_1V_1^{\vee}= (V_0^{\vee}V_1^{\vee})^2-
\bar{u}_0\bar{u}_1V_0^{\vee}V_1^{\vee}-\bar{u}_1V_1^{\vee}
-\bar{u}_0V_0^{\vee}-\bar{k}_0V_0^{\vee}V_1V_1^{\vee}.
\end{multline*}
here the first equality is (\ref{k0}), second is (\ref{VVVV})
third and fourth are (\ref{u0}), (\ref{u1});
\begin{equation*}
V_1^{-1}(V_0^{\vee})^{-1} V_0^{-1}V_1^{\vee}=V_1^{-1}V_1^\vee V_1
V_1^\vee=(V_1 V_1^\vee)^2-\bar{k}_1 V_1^\vee V_1 V_1^\vee,
\end{equation*}
here the first equation is (\ref{VVVV}), the second is (\ref{k1}).

Thus summing up the previous formulas and using  formulas
(\ref{u0}), (\ref{u1})  we get:
\begin{multline}\label{()X3}
((V^{\vee}_1)^{-1}V_0+V_0^{-1}V_1^{\vee})X_3= X_1^2+X_2^2-4
-\bar{u}^2_0-\bar{u}^2_1-\bar{u}_0\bar{u}_1 X_2+\\
\bar{k}_0((V_1^{\vee})^{-1}V_1^{-1}(V_0^{\vee})^{-1}-V_0^{\vee}
V_1V_1^{\vee})+
\bar{k}_1((V_1^\vee)^{-1}V_1^{-1}(V_1^{\vee})^{-1}-V_1^\vee V_1
V_1^\vee ).
\end{multline}
Expanding last two summands in the expression we get:
\begin{multline*}
(V_1^{\vee})^{-1}V_1^{-1}(V_0^{\vee})^{-1}-V_0^{\vee}
V_1V_1^{\vee}=V_1^{\vee}V_1(V_0^{\vee})^{-1}-V_0^{\vee}
(V_1)^{-1}(V_1^{\vee})^{-1}-\\
\bar{u}_1(V_1(V_0^\vee)^{-1}+V_0^\vee V_1^{-1})-\bar{k}_1(V_1^\vee
(V_0^\vee)^{-1}+V_0^\vee
(V_1^\vee)^{-1})+\bar{u}_1\bar{k}_1((V_0^\vee)^{-1}-V_0^\vee)=\\
(V_0^\vee)^{-1} (V_0)^{-1}(V_0^\vee)^{-1}-V_0^\vee V_0
V_0^\vee-\bar{u}_0\bar{u}_1\bar{k}_1-
\bar{u}_1(V_1(V_0^\vee)^{-1}+V_0^\vee V_1^{-1})-\\
\bar{k}_1(V_1^\vee (V_0^\vee)^{-1}+V_0^\vee (V_1^\vee)^{-1}),
\end{multline*}
here the first equality is (\ref{u1}), (\ref{k1}); the second
equality uses (\ref{VVVV}) twice and (\ref{u0}). Using (\ref{u0}),
(\ref{u1}), (\ref{k0}), (\ref{k1}) we can rewrite the last two
terms
\begin{gather*}
V_0^\vee V_1^{-1}+V_1(V_0^\vee)^{-1}=X_3+
\bar{u}_0(V_1^{-1}-V_1)=X_3-\bar{u}_0\bar{k}_1, \\
V_1^\vee(V_0^\vee)^{-1}+V_0^\vee(V_1^\vee)^{-1}=X_2+
\bar{u}_1((V_0^\vee)^{-1}-V_0^\vee)=X_2-\bar{u}_0\bar{u}_1,
\end{gather*}
and the first two terms:
\begin{multline*}
(V_0^\vee)^{-1} (V_0)^{-1}(V_0^\vee)^{-1}-V_0^\vee V_0 V_0^\vee=
V_0^\vee(V_0^{-1}-V_0)(V_0^\vee)^{-1}-\\ \bar{u}_0(V_0V_0^\vee+
(V_0^\vee)^{-1}V_0^{-1})=-\bar{k}_0-\bar{u}_0 X_1.
\end{multline*}
That is we get
\begin{multline}\label{1()X3}
(V_1^{\vee})^{-1}V_1^{-1}(V_0^{\vee})^{-1}-V_0^{\vee}
V_1V_1^{\vee}=-\bar{k}_0-\bar{u}_0-\bar{k}_1X_2-\bar{u}_1
X_3+\bar{u}_0\bar{u}_1\bar{k}_1.
\end{multline}

Now we expand the very last summand in (\ref{()X3}):
\begin{multline}\label{2()X3}
(V_1^\vee)^{-1}V_1^{-1}(V_1^{\vee})^{-1}-V_1^\vee V_1 V_1^\vee =
(V_1^\vee)^{-1}(V_1^{-1}-V_1)V_1^\vee-\\ \bar{u}_1(V_1V_1^\vee+
(V_1^\vee)^{-1}V_1^{-1})=-\bar{u}_1 X_1-\bar{k}_1.
\end{multline}
The first equality uses (\ref{u1}), the second uses (\ref{X1per})
and (\ref{k1}).

It easy to see that formulas (\ref{X1X2}), (\ref{()X3}),
(\ref{1()X3}), (\ref{2()X3}) imply the desired equation:
\begin{multline*}
X_1X_2X_3-X_1^2-X_2^2-X_3^2+(\bar{u}_0\bar{k}_0+\bar{k}_1\bar{u}_1)X_1+
(\bar{u}_1\bar{u}_0+\bar{k}_0\bar{k}_1)X_2+\\(\bar{k}_0\bar{u}_1+
\bar{k}_1\bar{u}_0)X_3+\bar{k}_0^2+\bar{k}_1^2+\bar{u}_0^2+\bar{u}_1^2-
\bar{k}_0\bar{k}_1\bar{u}_0\bar{u}_1+4=0.
\end{multline*}

\subsection{Proof of Proposition~\ref{Poisbr}}
 First
of all let us notice that the elements of the ring
$Z=\CC[C_{\underline{t}}]$ act by differentiation on $H$:
$$D_z(x)=[z_h,x]/h|_{h=0},$$ where $z_h$ is such that $z_0=z\in Z$
and $x\in H$. In particular, using formula (\ref{frml cntr}) we
get the action of $D_{X_1}$:
\begin{gather}
D_{X_1}(V_0)=(V_0V_0^{\vee}V_0^{-1}-V_0^{\vee})/2,\label{X1V0}\\
D_{X_1}(V_0^{\vee})=(V_0-(V_0^{\vee})^{-1}V_0V_0^{\vee})/2,\\
D_{X_1}(V_1)=(V_1^{\vee}-V_1^{-1}V_1^{\vee}V_1)/2,\\
D_{X_1}(V_1^{\vee})=(V_1^{\vee}V_1(V_1^{\vee})^{-1}-V_1)/2.\label{X1V1c}
\end{gather}
Let us show how to derive formula (\ref{X1V0}), the rest of the
formulas can be derived analogously: \begin{multline*}
[X_1,V_0]=V_1^\vee V_1 V_0+V_0V_0^\vee V_0- V_0 V_1^{\vee}
V_1-V_0V_0V_0^{\vee}=
q^{-1/2}(V_0^\vee)^{-1}V_0^{-1}V_0+\\V_0V_0^{\vee}V_0-q^{-1/2}V_0(V_0^\vee)^{-1}
V_0^{-1}-V_0^\vee-\bar{k}_0V_0V_0^\vee=-q^{-1/2}\bar{u}_0+(q^{-1/2}-1)V_0^\vee+\\
V_0(V_0^\vee-q^{-1/2}(V_0^\vee)^{-1})V_0^{-1}=q^{-1/4}(q^{1/4}-q^{-1/4})(V_0V_0^\vee
V_0^{-1}-V_0^\vee),
\end{multline*}
here we use (\ref{VVVV}) and (\ref{k0}) in the second equality; in
the third equality we use (\ref{u0}) and (\ref{k0}); in the fourth
equality we use (\ref{u0}) one more times.

From the formula for the differentiation we can derive
\begin{multline*}
2\{X_1,X_2\}_{def}=2D_{X_1}(X_2)=(V_1^\vee-V_1^{-1}V_1^\vee
V_1)V_0+V_1(V_0V_0^\vee
V_0^{-1}-V_0^\vee)+\\(V_0-(V_0^\vee)^{-1}V_0V_0^\vee)V_1^\vee+V_0^\vee
(V_1^\vee V_1(V_1^\vee)^{-1}-V_1)=(V_1^\vee V_0+V_0 V_1^\vee)-\\
(V_1V_0^\vee +V_0^\vee
V_1)+((V_1^\vee)^{-1}V_0^{-1}+V_0^{-1}(V_1^\vee)^{-1})-(V_1^{-1}
(V_0^\vee)^{-1}+(V_0^\vee)^{-1}V_1^{-1}).
\end{multline*}
Now let us expand each of the four terms in the RHS of the last
formula:
\begin{gather*}
V_1^\vee V_0+V_0V_1^\vee=(V_1^\vee)^{-1}V_0+V_0^{-1}V_1^\vee+
\bar{u}_1 V_0+\bar{k}_0 V_1^\vee,\\ V_1V_0^\vee +V_0^\vee
V_1V_1=X_3+\bar{u}_0V_1^{-1}+\bar{k}_1(V_0^\vee)^{-1}+\bar{u}_0\bar{k}_1,\\
(V_1^\vee)^{-1}V_0^{-1}+V_0^{-1}(V_1^\vee)^{-1}=(V_1^\vee)^{-1}V_0+
V_0^{-1}V_1^\vee-\bar{k}_0(V_1^\vee)^{-1}-\bar{u}_1 V_0^{-1},\\
V_1^{-1}
(V_0^\vee)^{-1}+(V_0^\vee)^{-1}V_1^{-1}=X_3-\bar{k}_1V_0^\vee-\bar{u}_0V_1+
\bar{u}_0\bar{k}_1),
\end{gather*}
here the first formula follows from (\ref{u0}) and (\ref{k1}); the
second  from (\ref{u1}), (\ref{k0}); the third  from (\ref{k0}),
(\ref{u1}); the last one from (\ref{u0}), (\ref{k1}). Thus we get
\begin{multline*}
2\{X_1,X_2\}_{def}=-2X_3+2((V_1^\vee)^{-1}V_0+V_0^{-1}V_1^\vee)+
\bar{u}_0(V_1-V_1^{-1})+\bar{k}_1(V_0^\vee-\\
(V_0^\vee)^{-1})+\bar{u}_1(V_0-V_0^{-1})+\bar{k}_0(V_1^\vee-(V_1^\vee)^{-1})-
2\bar{u}_0\bar{k}_1=-2X_3+2(X_1X_2-X_3+\\
\bar{u}_0\bar{k}_1)+2\bar{u}_1\bar{k}_0= 2\{X_1,X_2\},
\end{multline*}
here the second equality uses (\ref{k0}-\ref{u1}) and
(\ref{X1X2}).

As both brackets $\{\cdot,\cdot\}$ and $\{\cdot,\cdot\}_{def}$ are
algebraic and $X_1$, $X_2$ are local coordinates on the open part
of $C_{\underline{t}}$, the equality
$\{X_1,X_2\}=\{X_1,X_2\}_{def}$ implies the statement.

\end{document}